\documentclass[preprint,aos]{imsart}
\usepackage{amssymb}
\usepackage{amsmath}
\usepackage{amsthm}
\usepackage{graphicx}
\RequirePackage{multirow}

\RequirePackage[OT1]{fontenc}
\RequirePackage[numbers]{natbib}
\RequirePackage[colorlinks,citecolor=blue,urlcolor=blue]{hyperref}
\RequirePackage{hypernat}

\usepackage[left=1in,right=1.0in,top=1.0in,bottom=1.0in]{geometry}
\startlocaldefs

\newtheorem{conj}{Conjecture}
\newtheorem{thm}[conj]{Theorem}
\newtheorem{cor}[conj]{Corollary}
\newtheorem{prop}[conj]{Proposition}
\newtheorem{lemma}[conj]{Lemma}

\def\to{\rightarrow}

\def\EE{ {\rm I} \kern-.15em {\rm E} }
\def\PP{ {\rm I} \kern-.15em {\rm P} }

\def\wh{\widehat}
\def\eps{\varepsilon}

\def\X{\mbox{$\mathcal X$}}

\def\M{\mbox{$\mathcal M$}}
\def\N{\mbox{$\mathcal N$}}
\def\RR{\mathbb R}
\def\EE{\mathbb{E}}
\def\PP{\mathbb{P}}
\def\1{\mathbb{I}}

\def\argmin{\mathop{\rm arg\, min}}

\endlocaldefs

\begin{document}

\begin{frontmatter}
\title{ Optimal selection of reduced rank estimators of high-dimensional matrices}
\runtitle{High Dimensional Reduced Rank Estimator}

\begin{aug}
\author{\fnms{Florentina} \snm{Bunea}\thanksref{t1}\ead[label=e1]{flori@stat.fsu.edu}},
\author{\fnms{Yiyuan} \snm{She} \ead[label=e2]{yshe@stat.fsu.edu}} \and
\author{\fnms{Marten H.} \snm{Wegkamp}\thanksref{t1}\ead[label=e3]{wegkamp@stat.fsu.edu}}

\thankstext{t1}{The research of  Bunea and    Wegkamp was supported in part by NSF Grant DMS-1007444.}

\runauthor{F. Bunea, Y. She, and M.H. Wegkamp}

\affiliation{Florida State University}

\address{Department of Statistics\\
Florida State University\\
Tallahassee, FL 32306-4330\\
\printead{e1}\\
\phantom{E-mail:\ }\printead*{e2}\\
\phantom{E-mail:\ }\printead*{e3}
}

\end{aug}

\begin{abstract}
We introduce a new criterion, the Rank Selection Criterion (RSC), for selecting the optimal reduced rank estimator of the coefficient matrix in multivariate response regression models.  The corresponding RSC  estimator minimizes the Frobenius norm of the fit plus a regularization term  proportional to the number of parameters in the reduced rank model.

The rank of the RSC estimator provides a consistent estimator of  the rank of the coefficient matrix; in general the rank of our estimator is a consistent estimate of the {\it effective rank}, which we define to be the
number of singular values of the target matrix that are appropriately large. The consistency results are
valid not only in  the
classic asymptotic regime, when $n$, the number of responses, and $p$, the number of predictors, stay bounded, and $m$, the number of observations, grows, but also when either, or both, $n$ and $p$ grow, possibly much faster than  $m$.

We establish minimax optimal bounds
on the mean squared errors of our estimators.
Our finite sample   performance bounds for the RSC estimator show that it achieves the optimal balance between the approximation error and the penalty term.

Furthermore,  our procedure has very low computational complexity, linear in the number of candidate models, making it particularly appealing for large scale problems.
We contrast our estimator with  the nuclear norm penalized least squares (NNP) estimator, which has an inherently higher computational complexity than RSC, for multivariate regression models.  We show that  NNP has estimation properties similar to those of RSC, albeit under  stronger conditions. However,  it is not as parsimonious as RSC. We offer a simple correction of the NNP estimator  which leads to consistent rank estimation.

We verify and illustrate our theoretical findings via an extensive  simulation study.

\end{abstract}



\end{frontmatter}

\section{Introduction}
In  this paper we propose and analyze  dimension reduction-type estimators for   multivariate response regression models. Given $m$ observations of the responses $Y_i \in \RR^n$ and predictors $X_i \in \RR^p$, we assume that the matrices $Y=[Y_1,\ldots,Y_m]'$ and $X=[X_1,\ldots,X_m]'$ are related
via  an unknown $p\times n$ matrix of coefficients $A$, and write this as
\begin{eqnarray}
\label{mod}
 Y= XA +E,\end{eqnarray}
where $E$ is a random $m\times n$ matrix,  with independent  entries with mean zero and variance $\sigma^2$.

Standard least squares estimation in (\ref{mod}), under no constraints, is equivalent to regressing each response on the predictors separately. It  completely ignores the multivariate nature of the possibly correlated responses, see, for instance,  Izenman (2008) for a discussion of this phenomenon.  Estimators    restricted to have  rank equal to  a fixed number $k \leq n \wedge p$  were introduced to remedy this drawback. The history of such estimators dates back to the 1950's, and was initiated by Anderson (1951). Izenman (1975) introduced the term {\em reduced-rank
regression}  for this class of models and provided further study of the estimates. A number of important works followed, including Robinson (1973, 1974) and Rao (1978).
The monograph on reduced rank regression by Reinsel and Velu (1998) has  an excellent, comprehensive account of more recent developments and extensions of the model. All theoretical results to date for estimators of $A$ constrained to have rank equal to a given value $k$ are of asymptotic nature and are obtained  for fixed $p$, independent of the number of observations $m$.  Most of them are obtained in a likelihood framework, for Gaussian errors $E_{ij}$. Anderson (1999) relaxed this assumption and derived the asymptotic distribution of the estimate, when $p$ is fixed, the errors have two finite moments, and  the rank of $A$ is known.
Anderson (2002) continued this work by constructing  asymptotic tests for rank selection, valid only for small and fixed values of $p$.

The aim of our work is to develop a non-asymptotic class of  methods that yield reduced rank estimators of $A$ that are easy to compute, have rank determined adaptively from the data, and are valid for any values of $m, n$ and $p$, especially when the number of predictors $p$ is large.  The resulting estimators can then be used  to construct a possibly  much smaller number of new transformed predictors or can be used to construct the most important canonical variables based on the original $X$ and $Y$. We refer to Chapter 6 in Izenman (2008) for a  historical account of the latter.

We propose to  estimate $A$ by  minimizing  the sum of squares $\| Y- XB\|_F^2= \sum_i \sum_j \{ Y_{ij}- (XB)_{ij} \}^2$ plus a penalty $\mu r(B)$, proportional to the rank $r(B) $, over all matrices $B$.  It is immediate to see, using Pythagoras' theorem,  that  this is equivalent with computing
$ \min_{B} \left\{  \| PY - XB\|_F^2 + \mu r(B)  \right\}$
or $\min_k \left\{ \min_{B:\ r(B)=k}  \| PY - XB\|_F^2 +\mu k \right\}$,
with $P$ being  the projection matrix onto the column space of $X$.
In Section 2.1 we show that the minimizer $\wh k$ of the above expression is the number of  singular values $d_k( PY)$ of $PY$ that exceed  ${\mu}^{1/2}$. This observation reveals the prominent role of the tuning parameter $\mu$ in constructing $\wh k$. The final estimator $\wh A$  of the target matrix $A$ is the minimizer of $  \| PY - XB\|_F^2 $ over matrices $B$ of rank $\wh k$, and can be computed efficiently even for large $p$, using  the procedure that we describe in detail in Section 2.1 below.

The theoretical analysis of our proposed estimator $\wh A$ is  presented in Sections 2.2 -- 2.4.
The rank of $A$ may
not be the most appropriate measure of sparsity in multivariate regression models. For instance,
suppose that the rank of $A$ is 100, but only  three of its singular values are large and the remaining 97 are nearly zero.
This is an extreme example, and in general one needs an  objective method for  declaring
singular values as ``large" or ``small".  We introduce in Section 2.1  a slightly different
notion of sparsity, that of {\it effective rank}. The effective rank counts  the number of singular values of the signal $XA$ that are above a certain noise level.  The relevant notion of noise level  turns out to be
the largest singular value of $PE$. This is central to our results, and influences the choice of the tuning sequence $\mu$.  In Appendix C we prove that the expected value of the largest singular value of $PE$ is bounded by  $(q+n)^{1/2}$, where $q \leq m \wedge p$ is the rank of $X$. The effective noise level is at most $(m+n)^{1/2}$, for instance in the model $Y=A+E$,  but it can be substantially lower, of order $(q+n)^{1/2}$,  in  model (\ref{mod}).

In Section 2.2 we give tight conditions under which  $\wh k$, the rank of our proposed estimator $\wh A$, coincides with the effective rank. As an immediate corollary we show when  $\wh k$ equals the rank of $A$.  We give  finite sample performance bounds for $\|X\wh A - XA\|_F^2$
 in Section 2.3. These results
show  that $\wh A$ mimics the behavior of  reduced rank estimates based on the ideal effective rank, had this been known prior to estimation. If  $X$ has a restricted isometrity property, our estimate is minimax adaptive.
In the asymptotic setting, for $n+ (m\wedge p) \ge n+q \to\infty$, all our results hold with probability close to one, for  tuning parameter chosen proportionally to the square of the noise level.

We often particularize  our  main findings to the setting of  Gaussian $N(0,\sigma^2)$ errors $E_{ij}$ in order to obtain sharp, explicit numerical constants for the penalty term. To avoid technicalities, we assume that $\sigma^2$ is known in most cases, and we treat the case of unknown $\sigma^2$ in Section 2.4.

We contrast our estimator with the penalized least squares estimator $\widetilde{A}$ corresponding to a  penalty term $\tau \| B\|_1 $ proportional to the nuclear norm
$\| B\|_1=\sum_j d_j(B)$,  the sum of the singular values of $B$. This estimator has been studied by, among others,  Yuan et al. (2007) and  Lu et al (2010),  for model (\ref{mod}). Nuclear norm penalized estimators in
general models $y=\X(A)+\eps$ involving linear maps $\X$ have been studied by  Cand\`es and Plan (2010) and Negahban and Wainwright (2009).
A special case of this model  is the challenging
matrix completion problem, first investigated theoretically, in the noiseless case,  by
Cand\`es and Tao (2010). Rohde and Tsybakov (2010) studied a larger class of  penalized estimators, that  includes the nuclear norm estimator, in the general model $y=\X(A)+\eps$.

 In Section 3 we give bounds on $\|X\widetilde A - XA\|_F^2$
  that are similar in spirit to those from Section 2.  While the error bounds of the two estimators are comparable, albeit with cleaner results and milder conditions for our proposed estimator, there is one aspect in which the estimates differ in important ways.  The nuclear norm penalized estimator is far less parsimonious than the estimate obtained via our rank selection criterion. In Section 3, we offer a correction of the former estimate that yields a correct rank estimate.

Section 4 complements our theoretical results by an extensive simulation study that supports our theoretical findings and suggests strongly that the proposed estimator behaves very well in practice, in most situations is preferable to the nuclear norm penalized estimator and it is always much faster to compute.

Technical results and some intermediate proofs are presented in  Appendices A -- D.

\section{The Rank Selection Criterion}

\subsection{Methodology}

We propose to estimate $A$ by the penalized least squares  estimator
 \begin{eqnarray}
 \wh A = \argmin_{B} \{ \| Y- XB\|_F^2 + \mu r(B) \}.
 \end{eqnarray}
We denote its rank by $\wh k$.
The minimization is taken over all $p\times n$ matrices $B$. Here and in what follows
 $r(B)$ is the rank of $B$ and $\|C \|_F =  \left( \sum_i\sum_j C_{ij}^2 \right)^{1/2}$ denotes the  Frobenius norm for any generic matrix $C$. The choice of the  tuning parameter $\mu>0$ is discussed in Section 2.2.
  Since
\begin{equation} \label{unu}\min_{B}\left  \{ \| Y- XB\|_F^2 + \mu r(B) \right \} = \min_{k}\left \{ \min_{B, \ r(B) = k} \left   \{\| Y- XB\|_F^2 + \mu k \right \} \right \},\end{equation}
one needs to compute the  restricted rank estimators $\wh B_k$ that minimize $\| Y - XB\|_F^2$ over all matrices $B$ of rank $k$.
The following computationally efficient procedure for calculating each   $\wh B_k$ has been suggested by  Reinsel and Velu (1998).   Let $ M = X'X$ be the Gram matrix,  $M^{-}$ be its Moore-Penrose inverse and let  $P=X M^{-} X'$ be the projection matrix onto the column space of $X$.
\begin{enumerate}
\item Compute the eigenvectors $V=[v_1, v_2, \cdots, v_n]$, corresponding to the ordered eigenvalues arranged  from largest to smallest,   of  the symmetric matrix  $Y' P Y$.
\item    Compute the least squares estimator $\wh B= M^{-} X' Y$.\\  Construct  $W= \wh B V$
and $G = V^{\prime}$.\\ Form $ W_{  k} =  W[\,  ,1: k ]$ and
$ G_{  k} =  G[1:  k,\, ]$.
\item Compute the final estimator $\wh B_k = W_{  k}  G_{  k} $.
\end{enumerate}

In   step 2  above, $  W_{ k}$ denotes the matrix obtained from $ W$ by retaining all its rows and only its first $  k$ columns,  and  $G_{  k} $ is obtained from $G$ by retaining its first $ k$ rows and all its columns. \\

Our first result, Proposition \ref{ranksel} below,  characterizes the minimizer  $\wh k=r(\wh A)$ of (\ref{unu}) as the number of eigenvalues of the square matrix $Y'PY$ that exceed $\mu$ or, equivalently,  as the number of singular values of the matrix $PY$ that exceed ${\mu}^{1/2}$. The final estimator  of $A$  is then  $\wh A= \wh B_{\wh k}$.

Lemma \ref{lemma:GSVD} in Appendix B shows that the fitted matrix $X\wh A$ is equal to $\sum_{j\le \wh k} d_j u_j v_j'$ based on the singular value decomposition $UDV= \sum_{j} d_j u_j v_j'$ of the projection $PY$.

 \begin{prop}\label{ranksel}
Let $\lambda_1(Y'PY) \ge \lambda_2( Y'PY) \ge \cdots$ be the ordered eigenvalues of  $Y'PY$.
We have $\wh A= \wh B_{\wh k}$ with
\begin{equation}\label{khat} \wh k = \max\left\{ k:\ \lambda_k (Y'PY) \ge \mu \right\}.\end{equation}
\end{prop}
\begin{proof}

For $\wh B_k$ given above, and by the Pythagorean theorem, we have
\[  \| Y- X\wh B_k\|_F^2 = \| Y- PY\|_F^2 + \| PY - X\wh B_k\|_F^2, \]
and we observe that $X\wh B= PY$.
By Lemma \ref{lemma:GSVD} in Appendix B,
we have
\begin{eqnarray*}
\| X\wh B - X \wh B_k \|_F^2 = \sum_{j>k} d_j^2(X\wh B)
=  \sum_{j>k} d_j^2 (PY)= \sum_{j>k}\lambda_j( Y'PY),
\end{eqnarray*}
where $d_j(C)$ denotes the $j$-th largest singular value of a matrix $C$.
Then, the penalized least squares criterion reduces to
$$\| Y- PY\|_F^2 +  \left\{ \sum_{j > k} \lambda_j( Y'PY) +\mu k \right\}, $$
and we find that $\min_B \left\{ \| Y- XB\|_F^2 + \mu r(B) \right \}$ equals
 \begin{eqnarray*}
 \| Y- PY\|_F^2 - \mu n +  \min_k \sum_{j > k} \left\{   \lambda_j( Y'PY) - \mu \right\} .
\end{eqnarray*}
It is easy to see that $\sum_{j > k} \left\{   \lambda_j( Y'PY) - \mu \right\} $  is minimized by taking $k$ as the largest index $j$  for which $  \lambda_j( Y'PY) - \mu \ge 0$, since then the sum only consists of negative terms.
This concludes our proof.
\end{proof}

 {\sc Remark.}  The  two matrices $W_{  \wh k}$ and $G_{ \wh  k}$,  that yield
the final solution $\wh A = W_{  \wh k}G_{ \wh  k} $, have  the following properties:
(i) $G_ {\wh k} G_{ \wh  k} ^{\prime} $ is the identity matrix; and  (ii) $W_{\wh   k}^{\prime} M W_{ \wh  k}$ is a diagonal matrix.
Moreover,  the decomposition of $\wh A$    as a product of two matrices with  properties (i) and (ii) is unique, see, for instance, Theorem 2.2 in Reinsel and V\'elu (1998). As an immediate consequence, one can construct new orthogonal predictors as the columns of  $Z = XW_{\wh k}$. If $\wh k$ is much smaller than $p$, this can result in a significant dimension reduction of the predictors' space.
 \\

 \subsection{Consistent effective rank estimation}

In this section we study the properties of $\wh k=r(\wh A)$.
We will state simple conditions that guarantee that $\wh k$ equals $r =r(A)$ with high probability. First, we describe in Theorem \ref{rank_regression}
what $\wh k$ estimates and what quantities need to be controlled for consistent estimation. It turns out that  $\wh k$  estimates the number of the singular values of the signal $XA$ above  the threshold   $\mu^{1/2}$, for
 any value of the tuning parameter $\mu$.
The quality of estimation is controlled by the probability that   this threshold level exceeds the largest singular value $d_1(PE)$ of the {\em projected} noise matrix $PE$.
 We denote the $j$th singular value of a generic matrix $C$ by $d_j(C)$  and  we use the convention that the singular values are indexed in  decreasing order.

\begin{thm}\label{rank_regression}
Suppose that there exists an index  $s\le r$ such that  $$d_s(XA) > (1+\delta)\sqrt{\mu}\ \text{ and }
\ d_{s+1}(XA) < (1-\delta)\sqrt{\mu},$$  for some $\delta\in(0,1]$. Then
we have
\[ \PP\left\{ \wh k =  s\right\} \geq 1 -  \PP\left\{ d_1(PE) \ge \delta \sqrt{\mu} \right\} .\\
\]
\end{thm}
\begin{proof}
Using the characterization of $\wh k$ given in Proposition \ref{ranksel} we have
\begin{eqnarray*}
\wh k > s &\Longleftrightarrow & \sqrt \mu \leq d_{s+1} (PY)\\
\wh k < s &\Longleftrightarrow &\sqrt \mu \ge d_s (PY).
\end{eqnarray*}
Therefore
$ \PP\left\{ \wh k \ne s \right\} = \PP\left\{  \sqrt \mu \leq d_{s+1} (PY) \ \text{ or } \
\sqrt \mu \ge d_s (PY) \right\}.
$
Next, observe that $PY=XA+ PE$ and
$d_k(XA)<d_k(PY)+d_1(PE)$ for any $k$.
Hence $d_s(PY)\le \mu^{1/2}$ implies $d_1(PE) \ge d_s(XA)-\mu^{1/2}$, whereas
$d_{s+1}(PY) \ge \mu^{1/2}$ implies that $d_1(PE) \ge \mu^{1/2} - d_{s+1}(XA)$.
Consequently we have
\[ \PP\left\{ \wh k \ne s \right\} \le \PP\left\{ d_1(PE) \ge \min\left( \sqrt{\mu} - d_{s+1}({XA}), d_s({XA})-\sqrt\mu \right) \right\}.
\]
Invoke the conditions on $d_{s+1}({XA})$ and $d_{s}({XA})$ to complete the proof.
\end{proof}

Theorem \ref{rank_regression}  indicates that we can consistently estimate the index $s$
provided we use a large enough value for our tuning parameter  $\mu$ to guarantee that the probability of the event   $ \left\{ d_1(PE) \le \delta \mu^{1/2} \right\}$ approaches one.  We call $s$
the {\em effective rank} of $A$ relative to $\mu$, and denote it  by $r_e=r_e(\mu)$.

This is the appropriate notion of sparsity in the multivariate regression problem: we can only hope to recover those singular values of the signal $XA$ that are above the noise level $\EE[d_1(PE)]$. Their number, $r_{e}$, will be the target rank of the approximation of the mean response, and can be much smaller than $r = r(A)$.
We  regard the largest singular value $d_1(PE)$ as the relevant indicator  of the strength of the noise. Standard results on the largest singular value of Gaussian matrices show that $\EE [d_1(E)] \le \sigma (m^{1/2} + n^{1/2} )$ and similar bounds are available for subGaussian matrices, see, for instance, Rudelson and Vershynin (2010). Interestingly, the  expected value of the largest singular value $d_1(PE)$ of the projected noise matrix
is smaller: it is of order $(q+n)^{1/2}$ with $q=r(X)$. If $E$ has independent $N(0,\sigma^2)$ entries the following simple argument shows why this is the case.
\begin{lemma}\label{lemma:Gaussian}
Let $q=r(X)$ and assume that $E_{ij}$ are independent $N(0,\sigma^2)$ random variables.
Then \[ \EE \left[ d_1(PE) \right] \le \sigma \left( \sqrt{n} + \sqrt{q} \right)\]
and
\[ \PP \left\{ d_1(PE) \ge \EE [d_1(PE) ] + \sigma t \right\}\le \exp\left( - t^2/2 \right)
\] for all $t>0$.
\end{lemma}
\begin{proof}
Let $U\Lambda U'$ be the eigen-decomposition of $P$. Since $P$ is the projection matrix on the column space of $X$, only the first $q$ entries of $\Lambda$ on the diagonal equal to one, and all the remaining entries equal to zero. Then,
$d_1^2(PE)= \lambda_1( E' P E ) = d_1^2(\Lambda U'E)$. Since $E$ has independent  $N(0,\sigma^2)$ entries, the rotation $U'E$ has the same distribution as $E$. Hence
$\Lambda U'E$ can be written as a $q\times n$ matrix with Gaussian entries on top of  a $(m-q)\times n$ matrix of zeroes. Standard random matrix theory now states that $\EE [d_1( \Lambda U'E )] \le \sigma ({q}^{1/2} + {n}^{1/2} )$. The second claim of the lemma is a direct consequence of Borell's inequality, see, for instance, Van der Vaart and Wellner (1996),  after recognizing that $d_1(\Lambda U' E)$ is the supremum of a Gaussian process.
\end{proof}

In view of this result,  we take $\mu^{1/2} > \sigma ( n^{1/2} + q^{1/2})  $ as our measure of the noise level. The following corollary summarizes the discussion above and lists the main results of this section:  the proposed estimator based on the rank selection criterion (RSC) recovers consistently the effective rank $r_e$ and, in particular, the rank of $A$.

\begin{cor}\label{vier}
Assume that $E$ has independent $N(0,\sigma^2)$ entries.
For any $\theta>0$, set  $$\mu = (1+\theta)^2 \sigma^2 ( \sqrt{n}+\sqrt{q})^2 / \delta^2$$ with $\delta$ as in Theorem \ref{rank_regression}. Then
we have, for any $\theta>0$,
$$\PP\{ \wh k \ne  r_{e}(\mu)\} \le \exp\left( -\frac12 \theta^2 (n+q) \right) \to 0  \text{ as $q+n\to\infty$.}$$
 In particular, if  $d_r(XA)> 2\mu^{1/2} $ and $\mu^{1/2} = (1+\theta) \sigma (\sqrt{n} + \sqrt{q}) $, then
$$\PP\{ \wh k \ne  r\} \le \exp\left( -\frac12 \theta ^2 (n+q) \right) \to 0 \text{ as $q+n\to\infty$.}$$
\end{cor}

{\sc Remark.} Corollary \ref{vier} holds when  $q+n\to\infty$. If $q+n$ stays bounded, but    $m\to\infty$, the consistency results continue to hold when
 $q$ is replaced by $q\ln(m)$ in the expression of the tuning parameter $\mu$ given above. Lemma \ref{lemma:Gaussian} justifies this choice. The same remark applies to all theoretical  results in this paper.\\

{\sc Remark.}
A more involved argument is needed in order to establish the conclusion of Lemma  \ref{lemma:Gaussian} when  $E$ has independent
subGaussian entries.  We give this argument in
 Proposition \ref{een}  presented in Appendix C.  Proposition \ref{een}  shows, in particular,  that when  $ \EE[ \exp( t E_{ij}) ]\le \exp(t^2/\Gamma_E)$  for all $t>0,$ and for some $\Gamma_E<\infty$, we have
$$\PP\left\{ d_1^2(PE) \ge32 \Gamma_E (q+n) ({ \ln(5) + x })
 \right\} \le 2 \exp\{-x(q+n)\},$$
for  all $x>0$. The conclusion of Corollary \ref{vier}  then holds for $\mu = C_0 \Gamma_E (n+q)$ with $C_0$ large enough.  Moreover, all oracle inequalities presented in the next sections remain valid  for this choice of the tuning parameter, if $E$ has independent
 subGaussian entries.
\\

 \subsection{Errors bounds for the RSC estimator}

In this section we study the performance of $\wh A$ by obtaining  bounds for  $\|X\wh A - XA\|_F^2$.  First we derive a bound for the fit $\| X\wh B_k - XA\|_F^2 $, based
on the restricted rank estimator $\wh B_k$, for each value of $k$.

 \begin{thm}\label{eachk}
Set $c(\theta)= 1+2/\theta$.
For any $\theta>0$, we have
\begin{eqnarray*}
\| X\wh B_k - XA\|_F^2 \le   \left\{ c^2(\theta)   \sum_{j>k} d_j^2(XA)
+  2(1+\theta) c(\theta)    k d_1^2(PE) \right\}
\end{eqnarray*}
with probability one.
\end{thm}
\begin{proof}
By the definition of $\wh B_k$,
\[ \| Y- X\wh B_k\|_F^2   \le \| Y- XB\|_F^2  \]
for all $p\times n$  matrices $B$ of rank $k$. Working out the squares we obtain
\begin{eqnarray*}
\| X\wh B_k - XA \|_F^2 & \le &
 \| XB-XA\|_F^2   + 2<E,X\wh A- XB>_F \\
 &=&   \| XB-XA\|_F^2   + 2<PE,X\wh A- XB>_F
\end{eqnarray*}
with
 \[ <C, D>_F= tr( C'D) =tr(D'C) =\sum_i \sum _j C_{ij} D_{ij}, \]
for generic $m \times n$ matrices $C$ and $D$.
The inner product $<C,D>_F$, operator norm $\|C\|_2= d_1(C)$ and nuclear norm
  $ \|D\|_1 = \sum_{j}d_j(D)$
are related via the inequality
  $<C,D> _F \le \| C\|_2 \| D\|_1$.
As a consequence we find
\begin{eqnarray*}
< PE, X\wh B_k- XB>_F &\le& d_1(PE) \|X\wh B_k- XB\|_1 \\
&\le&d_1(PE) \sqrt{ 2k  } \| X\wh B_k-XB\|_F\\
&\le& d_1(PE) \sqrt{2k}\{ \| X\wh B_k-XA\|_F + \| XB -XA\|_F\}.
\end{eqnarray*}
Using  the inequality $2xy \le x^2/a + ay^2$  with $a>0$ twice,
we obtain that $\| X\wh B_k - XA \|_F^2$ is bounded above by
\begin{eqnarray*}
\frac{1+b}{b}  \| XB-XA\|_F^2 +
\frac{1}{a}  \|X\wh B_k -XA \|_F^2+  (a+b)  (2k) d_1^2(PE).
\end{eqnarray*}
Hence  we obtain, for any $a,b>0$, the inequality
\begin{eqnarray*}
\| X\wh B_k - XA \|_F^2 & \le &
\frac{a}{a-1} \left\{  \frac{1+b}{b}  \| XB- XA \|_F^2 +   2(a+b) k d_1^2(PE) \right\}.
\end{eqnarray*}
Lemma \ref{lemma:GSVD} in the Appendix B states that the minimum of
$\| XA- XB\|_F^2 $ over all matrices $B$ of rank $k$ is achieved for the GSVD of $A$ and the minimum equals
$ \sum_{j>k} d_j^2(XA) $.
The claim follows after choosing
$a=(2+\theta)/2$ and $b=\theta/2$.
\end{proof}

\begin{cor}\label{eachk2}
Assume that $E$ has independent $N(0,\sigma^2)$ entries.
Set $c(\theta)= 1+2/\theta$. Then,  for any $\theta,\xi>0$, the inequality
\begin{eqnarray*}
&& \| X\wh B_k - XA\|_F^2  \\ &&\le   \left\{ c^2(\theta)  \sum_{j>k} d_j^2(XA)
+ 2 c(\theta) (1+\theta) (1+\xi)^2      \sigma^2 k  (n+q) \right\}
\end{eqnarray*} holds with probability $1-\exp(-\xi^2 (n+q) /2)$.
In addition,
\begin{eqnarray*}
\EE \left[  \| X\wh B_k - XA\|_F^2  \right]  \lesssim   \sum_{j>k} d_j^2(XA)
+       \sigma^2 k  (n+q).
\end{eqnarray*} The symbol $\lesssim$ means that the inequality holds up to  multiplicative numerical constants.
\end{cor}
\begin{proof}  Set $t= (1+\xi)^2 \sigma^2 (\sqrt{n}+\sqrt{q})^2$ for some $\xi>0$.
From Lemma \ref{lemma:Gaussian}, it follows that
\[ \PP\{ d_1^2(PE) \ge t\} = \PP\{ d_1(PE) \ge (1+\xi) \sigma (\sqrt{n}+\sqrt{q}) \} \le \exp(-\xi^2 (n+q)/2).\]
The first claim follows now from this bound and  Theorem \ref{eachk}.
From Lemma \ref{basis1}, it follows that $\EE[ d_1^2(PE) ] \le \nu^2 + \nu\sqrt{2\pi} + 2 $ for $\nu= \EE[ d_1(PE)]\le \sigma(\sqrt{n}+\sqrt{q})$. This proves the second claim.
\end{proof}

Theorem \ref{eachk}  bounds the error $\| X\wh B_k - XA\|_F^2$ by an approximation error, $\sum_{j>k} d_j^2(XA)$, and a stochastic term, $k d_1^2(PE)$, with probability one. The approximation error is decreasing in $k$ and vanishes for $k> r(XA)$.

The stochastic term  increases in $k$ and can be bounded by a constant times $ k (n+q)$ with overwhelming probability and in expectation, for Gaussian errors, by Corollary \ref{eachk2} above.  More generally, the same bound (up to constants) can be proved for  subGaussian errors.
 Indeed, for $C_0$ large enough, Proposition \ref{een} in Appendix  C, states that $\PP\{ d_1^2(PE) \le C_0 (n+q) \} \le2 \exp\{ - (n+q)\}$.

We observe that $k(n+q)$ is essentially the number of free parameters of the restricted rank problem. Indeed, our parameter space consists of all  $p\times n$ matrices $B$ of rank $k$ and each $XB$ matrix has
$k(n+q-k)$ free parameters. Hence we can interpret the bound in Corollary \ref{eachk2}  above as the squared bias plus the dimension of the parameter space.

Remark{\it (ii)}, following Corollary \ref{oracle1a} below,  shows that $k(n+q)$ is also the minimax lower bound for $\| X\wh B_k - XA\|_F^2$, if the smallest eigenvalue of $X'X$ is larger than a strictly positive constant. This means that  $X\wh B_k$   is a minimax estimator under this assumption.
\\

 We now turn to  the penalized estimator $\wh A$ and show that it  achieves the best (squared) bias-variance trade-off among all rank restricted estimators $\wh B_k$ for the appropriate choice of the tuning parameter $\mu$ in the penalty $\text{pen}(B)= \mu r(B)$.

\begin{thm}\label{oracle1}
We have, for any  $\theta>0$,  on the event  $ (1+\theta) d_1^2(PE)\le \mu$,
\begin{eqnarray}
  \| X\wh A - XA \|_F^2 \le  c^2(\theta)  \| XB-XA\|_F^2 + 2c(\theta) \mu k,
  \end{eqnarray} for any $p\times n$ matrix $B$.
In particular,
we have, for  $\mu\ge  (1+\theta) d_1^2(PE)$
\begin{eqnarray}\label{long}
\| X\wh A - XA\|_F^2 \le  \min_k  \left\{  c^2(\theta) \sum_{j>k} d_j^2(XA)
+  2 c(\theta) \mu k\right\}
\end{eqnarray}
and
\begin{eqnarray}\label{ranko}
\| X\wh A - XA\|_F^2 \le   2c(\theta) \mu   r.
\end{eqnarray}
\end{thm}
\begin{proof}
By definition of $\wh A$,
\[ \| Y- X\wh A\|_F^2 + \mu r(\wh A) \le \| Y- XB\|_F^2 + \mu r(B)\]
for all $p\times n$  matrices $B$. Working out the squares we obtain
\begin{eqnarray*}
&& \| X\wh A - XA \|_F^2 \\ && \le
 \| XB-XA\|_F^2 + 2\mu r(B) + 2<E,X\wh A- XB>_F-\mu r(\wh A)-\mu r(B)\\
 && =
 \| XB-XA\|_F^2 + 2\mu r(B) + 2<PE,X\wh A- XB>_F-\mu r(\wh A)-\mu r(B).
\end{eqnarray*}
Next we observe that
\begin{eqnarray*}
&& < PE, X\wh A- XB>_F \\
&&\le d_1(PE) \|X\wh A- XB\|_1 \\
&&\le d_1(PE) \{ r(X\wh A)+r(XB) \}^{1/2} \| X\wh A-XB\|_F\\
&&\le d_1(PE) \{ r(\wh A)+r(B) \}^{1/2}\{ \| X\wh A-XA\|_F + \| XB -XA\|_F\}.
\end{eqnarray*}
Consequently, using  the inequality $2xy \le x^2/a + ay^2$  twice,
we obtain, for any $a>0$ and $b>0$,
\begin{eqnarray*}
\| X\wh A - XA \|_F^2 & \le &
\| XB-XA\|_F^2 +
\frac{1}{a} \|X\wh A-XA \|_F^2 +\frac{1}{b}  \| XB-XA\|_F^2 +\\
&& 2 \mu r(B)  + (a+b)\{ r(\wh A)  + r( B)\} d_1^2(PE)
-\mu\{ r(\wh A)+r(B)\}.
\end{eqnarray*}
Hence, if $(a+b)d_1^2(PE)
-\mu\le 0$, we obtain
\begin{eqnarray*}
\| X\wh A - XA \|_F^2 & \le &
\frac{a}{a-1}\left\{ \frac{1+b}{b} \| XB- XA \|_F^2 + 2 \mu r(B)\right\},
\end{eqnarray*} for any $a>1$ and $b>0$.
Lemma \ref{lemma:GSVD} in  Appendix B evaluates the   minimum of
$\| XA- XB\|_F^2 $ over all matrices $B$ of rank $k$ and shows that it  equals
$ \sum_{j>k} d_j^2(XA) $.
We conclude our proof by choosing $a=1+\theta/2$ and $b=\theta/2$.
\end{proof}

{\sc Remark.}
The first two parts of the theorem  show that $\wh A$ achieves the best (squared) bias-variance trade-off among all reduced rank estimators $\wh B_k$ if $\mu>d_1^2(PE)$. Moreover, the index $k$ which minimizes  $\sum_{j>k} \{d_j^2(XA) +\mu k\}$ essentially coincides with the effective rank $r_e=r_e(\mu)$ defined in the previous section. Therefore, the fit of the selected estimator $X\wh A$ is comparable with that of the estimator $X\wh B_k$ with rank $k=r_e$. Since the ideal $r_e$  depends on the unknown matrix $A$, this ideal estimator cannot be computed. Although  our estimator $\wh A$ is constructed independently of $r_e$, it  mimics the  behavior of the ideal estimator $\wh B_{r_e}$ and we say that the bound on $\|X\wh A - XA\|_F^2$ adapts to $r_e \leq r$.

The last part of  our result is a particular case of the second part, but it is perhaps easier to interpret.
Taking the index $k$ equal to the rank $r$, the bias term disappears and the bound reduces to $r  d_1^2(PE)$ up to constants.
 This shows clearly the important role played by $r$ in  the estimation accuracy: the smaller the rank of $A$, the smaller the estimation error.  \\

For Gaussian errors, we have the following precise bounds.

\begin{cor}\label{oracle1a}
Assume that $E$ has independent $N(0,\sigma^2)$ entries. Set
\[ pen(B)= (1+\theta)(1+\xi) ^2 (\sqrt{n}+\sqrt{q})^2 \sigma^2 r(B)\]
with $\theta, \xi >0$ arbitrary. Let $c(\theta)=1+2 / \theta$.
Then, we have
\begin{eqnarray*}
\PP\left[ \| X\wh A - XA \|_F^2 \le  \min_{1\le k\le \min(n,p)} \left\{  c^2({\theta}) \sum_{j>k}d_j^2(XA) + 2c(\theta)\mu k\right\} \right]\\
 \ge 1- \exp\left\{ -\frac{ \xi^2 (n+q)}{2} \right\}
\end{eqnarray*}
and
\begin{eqnarray*}
&&\EE\left[
\| X\wh A - XA \|_F^2 \right]  \\ &&\le
\min_{1\le k\le \min(n,p)}  \left[ c^2(\theta) \sum_{j>k} d_j^2(XA)
 +   2(1-\theta)c(\theta) (1+\xi)^2 \sigma^2 (\sqrt{n}+\sqrt{q})^2 k \right] \\
&&\qquad +4(1+\theta)c(\theta) \min(n,p) \sigma^2 (1+\xi^{-1})  \exp(-\xi^2 (n+q)).
\end{eqnarray*}
\end{cor}

\begin{proof}
Recall from the proof of Theorem \ref{oracle1} that
\begin{eqnarray*}
\| X\wh A - XA \|_F^2 & \le & \frac{2+\theta}{\theta} \left\{  \frac{2+\theta}{\theta}
\| XB-XA\|_F^2 + 2 \text{pen} (B)  + R \right\}
\end{eqnarray*}
with $R$ defined by
\begin{eqnarray*}
R&=&
(1+\theta)\{ r(\wh A)  + r( B)\} d_1^2(PE)
- \text{pen} (\wh A) -\text{pen}(B) \\
&\le& 2 (1+\theta) \max_{1\le k\le \min(n,p) } k \left\{  d_1^2(PE)
- (1+\xi)^2 (\sqrt{n}+\sqrt{q}) ^2 \sigma^2  \right\}.
\end{eqnarray*}
For $\tilde E= E /\sigma$, a matrix of independent $N(0,1)$ entries, we have
\begin{eqnarray*}
R &\le& 2 (1+\theta) \sigma^2 \max_{1\le k\le \min(n,p) } k \left\{  d_1^2(P\tilde E)
- (1+\xi)^2 (\sqrt{n}+\sqrt{q}) ^2 \right\}\\
&\le& 2 \min(n,p) (1+\theta) \sigma^2 \left(  d_1^2(P\tilde E)
- (1+\xi)^2 (\sqrt{n}+\sqrt{q}) ^2 \right)_{+}.
\end{eqnarray*}
Apply Lemma \ref{basis1} in Appendix D to deduce that
\[ \EE [R] \le 4\min(n,p) \frac{1+\xi}{\xi} (1+\theta) \sigma^2 \exp(-\xi^2(n+q)/2).
\]
The conclusion  follows immediately.
\end{proof}

 {\sc Remarks.}
(i) We note that  for $n+q$ large,
\begin{eqnarray*}
&&\EE\left[
\| X\wh A - XA \|_F^2 \right] \lesssim
\min_{1\le k\le \min(n,p)}  \left[  \sum_{j>k} d_j^2(XA)
 +    \sigma^2 (\sqrt{n}+\sqrt{q})^2 k \right]
\end{eqnarray*}
as the remainder term in the bound  of $\EE\left[
\| X\wh A - XA \|_F^2 \right] $ in Corollary \ref{oracle1a}  converges exponentially fast in $n+q$, to zero.

(ii)  Assuming that  $E$ has independent $N(0,\sigma^2)$ entries,
the RSC estimator  corresponding to  the penalty $\text{pen}(B)=C \sigma^2 (n^{1/2} + q^{1/2})^2 r(B)$, for any
$C>1$, is minimax adaptive, for matrices $X$ having a restricted isometry property (RIP),
of the type introduced and discussed in Cand\` es and Plan (2010) and  Rohde and Tsybakov (2010). The RIP implies that $\|XA\|_{F}^{2} \geq \rho \|A\|_F^2$, for all matrices $A$ of rank at most $r$ and  for some   constant $ 0< \rho < 1$.  For fixed design matrices $X$,  this is equivalent with assuming that  the smallest eigenvalue $\lambda_p(M)$ of the $p\times p$ Gram  matrix $M = X^{\prime}X$ is  larger than $\rho$.
To establish the minimax lower bound for the mean squared error $\| X\wh A - XA\|^2_F$, notice first that  our model (\ref{mod}) can be rewritten  as  $y_i = \text{trace}(Z_i^{\prime}A) + \varepsilon_i$,
with $1 \leq i \leq mn$,  via the mapping $(a, b) \rightarrow i = a + (b - 1)n$, where
$ 1 \leq a \leq m$, $1 \leq b \leq n$, $y_i =: Y_{ab} \in \RR$ and $Z_i =: X_a^{\prime}e_b \in M_{p \times n}$.
Here $X_a \in \RR^p$ denotes the  $a$-th row of $X$, $e_b$ is the row vector in $\RR^n$ having the $b$-th component equal to 1 and the rest equal to zero, and $M_{p\times n}$ is the space of all $p\times n$ matrices. Then, under RIP, the  lower bound follows directly from Theorem 5 in Rohde and Tsybakov (2010); see also Theorem 2.5 in Cand\`es and Plan (2010) for minimax lower bounds on $\|\wh A - A \|_F^2$.

(iii)  The same type of upper bound as the one of Corollary \ref{oracle1a} can be proved if the entries of $E$ are subGaussian: take $\text{pen}(B)=C  (n+q)  r(B)$ for some $C$ large enough, and invoke Proposition \ref{een} in Appendix C.

(iv)
Although the error bounds of $\|X\wh A-XA\|_F$ are guaranteed for all $X$ and $A$, the analysis of the estimation performance  of $\wh A$  depends on $X$. If $\lambda_p(M) \geq \rho > 0$, for some constant $\rho$,   then, provided $\mu>  (1+\theta) d_1^2(PE)$ with $\theta>0$ arbitrary,
\begin{eqnarray*}
 \| \wh A - A\|_F^2 \le \frac{c(\theta)}{\lambda_{p}(M) }
\min_{k\le r} \left[c(\theta)  \sum_{j>k} d_j^2(XA)
+2\mu k \right]  \end{eqnarray*}
follows from Theorem \ref{oracle1}.

(v)
Our results are  slightly more general than stated. In fact, our analysis does  not require that the postulated multivariate linear model $Y=XA+E$ holds exactly.
We denote the expected value of $Y$ by $\Theta$  and write $Y=\Theta +E$. We denote  the projection of $\Theta$ onto the column space of $X$ by $XA$, that is, $P \Theta = XA$.
Because minimizing  $\| Y- XB\|_F^2  + \mu r(B)$  is equivalent with minimizing
$\| PY- XB\|_F^2  + \mu r(B)$ by Pythagoras' theorem, our least squares procedure estimates $XA$, the mean of $PY$.
The statements of
Theorems \ref{rank_regression} and \ref{oracle1}  remain unchanged, except that   $XA$ is  the mean of the projection $PY$ of $Y$, not the mean of $Y$ itself.\\

\subsection{A data adaptive penalty term}
\label{subsec:unknownvar}

In this section we construct a data adaptive penalty term that employs the unbiased estimator  \[ S^2 = \| Y- PY\|_F^2 / (mn-qn)\]
of $\sigma^2$.  Set, for any $\theta>0$, $\xi>0$ and $0<\delta<1$,
\[ \text{pen}(B) = \frac{(1+\theta)}{1-\delta}(1+\xi) ^2 (\sqrt{n}+\sqrt{q})^2 S^2 r(B). \]
 Notice that the estimator $S^2$ requires that $n(m-q)$ be large, which  holds whenever $m>>q$ or $m-q\ge1$ and $n$ is large. The challenging case $m=q<<p$ is left for future research.

\begin{thm}
Assume that $E$ is an $m\times n$ matrix with independent $N(0,\sigma^2)$ entries. Using the penalty given above
we have, for
 $c(\theta)= 1+ 2 / \theta$,
\begin{eqnarray*}
&&\EE\left[
\| X\wh A - XA \|_F^2 \right] \\
&&\le
\min_{1\le k\le \min(n,p)}  \left[ c^2(\theta) \sum_{j>k} d_j^2(XA)
 +   2(1+\theta) c(\theta) (1+\xi)^2 \sigma^2 (\sqrt{n}+\sqrt{q})^2 k \right] \\
 &&\quad+ 4(1+\theta) c(\theta) \min(n,p) \sigma^2(1+\xi^{-1}) \exp\left (-\frac{\xi^2({n}+{q})}{2} \right)\\
 &&\quad +  4(1+\theta) c(\theta) \min(n,p) \sigma^2\left( 2+  (\sqrt{n}+\sqrt{q})^2 + (\sqrt{n}+\sqrt{q}) \sqrt{2\pi}  \right) \times\\ &&\qquad\times \exp\left\{ -\frac{ \delta^2 n(m-q)}{4(1+\delta)}  \right\}.
 \end{eqnarray*}
\end{thm}
\begin{proof} Set $\tilde E=\sigma^{-1} E$.
We have, for any $p\times n$ matrix $B$
\begin{eqnarray*}
&&  \| X\wh A - XA \|_F^2\,  \le\,  \frac{2+\theta}{\theta} \left[ \frac{2+\theta}{\theta}
\| XB-XA\|_F^2 + 2 \text{pen}(B) \right] \\ &&\qquad +
2 \frac{2+\theta}{\theta}(1+\theta) \sigma^2  \max_{1\le k\le \min(n,p) } k \left\{  d_1^2(P \tilde E) -
\frac{ (1+\xi)^2 (\sqrt{n}+\sqrt{q})^2  S^2}{(1-\delta)\sigma^2} \right\}.
\end{eqnarray*}
It remains to bound the expected value of
\begin{eqnarray*}
&&\max_{k\le \min(n,p) } k \left\{  d_1^2(P \tilde E) -  (1+\xi)^2 (\sqrt{n}+\sqrt{q})^2\frac{S^2}{(1-\delta)\sigma^2} \right\}
\\ &&\le \min(n,p)  \left(  d_1^2(P \tilde E) -  (1+\xi)^2 (\sqrt{n}+\sqrt{q})^2 \frac{S^2}{(1-\delta)\sigma^2} \right)_{+}.
\end{eqnarray*}
We split the expectation into two parts: $S^2\ge (1-\delta)\sigma^2$ and its complement. We observe first that
\begin{eqnarray*}
&&\EE\left[ \left(  d_1^2(P \tilde E) -  (1+\xi)^2  (\sqrt{n}+\sqrt{q})^2\frac{S^2}{(1-\delta)\sigma^2} \right)_{+}1_{\{ S^2\ge (1-\delta)\sigma^2\}}
\right]
\\
&&\le
\EE\left[ \left(  d_1^2(P \tilde E) -  (1+\xi)^2 (\sqrt{n}+\sqrt{q})^2  \right)_{+}
\right]
\\
&&\le2 (1+{\xi}^{-1}) \min (n,p) \exp(-\xi (\sqrt{n}+\sqrt{q})/2),
\end{eqnarray*}
using Lemma \ref{basis1} for the last inequality.
Next, we observe that
\begin{eqnarray*}
&&\EE\left[ \left(  d_1^2(P \tilde E) -  (1+\xi)^2 (\sqrt{n}+\sqrt{q})^2 \frac{S^2}{(1-\delta)\sigma^2} \right)_{+}1_{\left\{ S^2\le (1-\delta)\sigma^2\right\}}
\right]
\\
&&\le \EE\left[   d_1^2(P \tilde E)1_{\left
\{ S^2\le (1-\delta)\sigma^2 \right\}}
\right]\,=\, \EE\left[   d_1^2(P \tilde E)
1_{\left\{ \| (I-P)\tilde E\|_F^2 \le (1-\delta)(nm-nq)\right\}}
\right].
\end{eqnarray*}
Since $P\tilde E$ and $(I-P)\tilde E$ are independent, and $\| (I-P)\tilde E\|_F^2$ has a $\chi^2_{nm-nq}$ distribution, we find
\begin{eqnarray*}
&&\EE\left[ \left(  d_1^2(P \tilde E) -  (1+\xi)^2 (\sqrt{n}+\sqrt{q})^2 \frac{S^2}{(1-\delta)\sigma^2} \right)_{+}
1_{\left\{ S^2\le (1-\delta)\sigma^2 \right\}}
\right]
\\
&&\le  \EE\left[  d_1^2(P \tilde E) \right] \PP
\left\{ \| (I-P)\tilde E\|_F^2 \le (1-\delta)(nm-nq)
\right\}\\
&&\le \left(   (\sqrt{n}+\sqrt{q})^2 +   \sqrt{2\pi}   (\sqrt{n}+\sqrt{q}) +2  \right) \exp\left\{ -\frac{ \delta^2}{4(1+\delta)} n(m-q) \right\},
\end{eqnarray*}
using Lemmas \ref{basis1} and \ref{chisq} in Appendix D for the last inequality. This proves the result.
\end{proof}

{\sc Remark.}
We see that for large values of $n+q$  and $n(m-q)$,
\begin{eqnarray*}
&&\EE\left[
\| X\wh A - XA \|_F^2 \right] \lesssim
\min_{1\le k\le \min(n,p)}  \left[   \sum_{j>k} d_j^2(XA)
 +    \sigma^2 (\sqrt{n}+\sqrt{q})^2 k \right],
\end{eqnarray*}
as the additional terms in the theorem above decrease exponentially fast in $n+q$  and $n(m-q)$.
This bound is similar to the one in Corollary \ref{oracle1a}, obtained  for the RSC estimator corresponding to the penalty term that employs  the theoretical  value of $\sigma^2$.\\

 \section{Comparison with nuclear norm  penalized estimators}
  In this section we compare our RSC estimator $\wh A$ with the  alternative estimator $\widetilde A$ that minimizes
\[ \| Y- XB\|_F^2 + 2 \tau \| B\|_1\]
over all $p\times n$ matrices $B$.

\begin{thm}\label{thm:NNPoracle1}
On the event $d_1(X' E) \le \tau$, we have, for any $B$,
\begin{eqnarray*}\|X\widetilde A - XA \|_F^2 \le \| XB-XA\|_F^2 + 4 \tau \| B\|_1.
\end{eqnarray*}
\end{thm}
\begin{proof}
By the definition of
$\widetilde A$,
\[ \| Y- X\widetilde  A\|_F^2 + 2\tau \|\widetilde  A\|_1  \le \| Y- XB\|_F^2 + 2\tau \|B \|_1 \]
for all $m\times n$  matrices $B$. Working out the squares we obtain
\begin{eqnarray*}\label{square1}
\| \widetilde XA -XA \|_F^2 & \le &
\| XB-XA\|_F^2 +  2\tau \|B \|_1 + 2<X'E,\widetilde  A-B>_F- 2\tau \|\widetilde  A\|_1
\end{eqnarray*}
Since
\begin{eqnarray*}\label{inner}
<X' E,\widetilde A-B>_F \le \| X' E\|_2 \| \widetilde A-B\|_1 \le {\tau}\|\widetilde{A} - B\|_1  \nonumber \end{eqnarray*}
on the event $d_1(X'E)\le  {\tau}$,
we obtain
the claim
using the triangle inequality.
\end{proof}

We see that $\widetilde A$ balances
the bias term $\| XA-XB\|_F^2$ with the penalty term  $\tau \| B\|_1$, provided $\tau > d_1(X'E)$. Since $X'E = X' PE+X'(I-P)E= X'PE$, we have
$d_1(X'E)\le d_1(X) d_1(PE)$.  We  immediately obtain the following corollary using the results for $d_1(PE)$ of Lemma \ref{lemma:Gaussian}.

\begin{cor} Assume that $E$ has independent $N(0,\sigma^2)$ entries. For
$$\tau = (1+\theta) d_1(X) \sigma (\sqrt{n}+\sqrt{q}) $$
with $\theta>0$ arbitrary, we have
\begin{eqnarray*}
\PP \left\{ \|X\widetilde A - XA \|_F^2 \le \| XB-XA\|_F^2 + 4 \tau \| B\|_1 \right\}
\ge 1-\exp\left\{ -\frac12 \theta^2 (n+q) \right\}.
\end{eqnarray*}
\end{cor}

The same result, up to constants, can be obtained if the errors $E_{ij}$ are subGaussian, if we  replace $\sigma$ in the choice of $\tau$ above by a suitably large constant $C$.
The proof of this generalization uses Proposition \ref{een} in Appendix C in lieu of Lemma \ref{lemma:Gaussian}. The same remark applies for all the results in this section.\\

The next result obtains an oracle inequality for $\widetilde A$  that resembles the oracle inequality for the RSC estimator $\widehat A$ in Theorem \ref{oracle1}. We stress the fact that Theorem \ref{Th-again} below requires that $\lambda_{p}(X'X) >0$;   this was  not required for the derivation of the oracle bound on $\| X\widehat  A- XA\|_F^2$  in Theorem \ref{oracle1}, which holds for all $X$.  We denote the condition number of $M=X'X$ by $c_0(M)=\lambda_1(M)/ \lambda_p(M)$.

\begin{thm}\label{Th-again}
 Assume that $E$ has independent $N(0,\sigma^2)$ entries. For
$$\tau = (1+\theta) d_1(X) \sigma (\sqrt{n}+\sqrt{q}) $$
with $\theta>0$ arbitrary, we have
\[ \| X\widetilde A- XA\|_F^2 \lesssim  \min_{k \le r} \left(  \sum_{j=k+1}^r d_j^2 (XA) + c_0(M)  k \sigma^2 (n+q)   \right)  .\]
Furthermore,
 \[ \| \widetilde A- A\|_F^2 \lesssim   c_0(M)  \sum_{j=k+1}^r d_j^2 (A) +  \frac{c_0(M) }{\lambda_{p}(M) } k \sigma^2 (n+ q)
     .\]
 Both inequalities hold with probability at least $1-\exp\left(-\theta^2 (n+q) /2 \right)$. The symbol $\lesssim$ means that the inequality holds up to  multiplicative numerical constants (depending on $\theta$).
\end{thm}
To keep the paper self contained, we give a simple proof of this result in Appendix A. Similar results for the NNP estimator of $A$ in the general model $y = \X(A)+ \varepsilon$, where  $\X$ is a random linear map,  have been obtained by
 Negahban and Wainwright  (2009) and Cand\`es and Plan (2010), each under different sets of assumptions on $\X$. We refer to Rohde and Tsybakov (2010)
 for more general results on Schatten norm penalized estimators of $A$ in the model $y = \X(A)+ \varepsilon$, and a very thorough discussion on the assumptions on $\X$ under which these results  hold. \\

\begin{figure}[htp!]
\begin{center}
\includegraphics[width=5cm, height=5cm ]{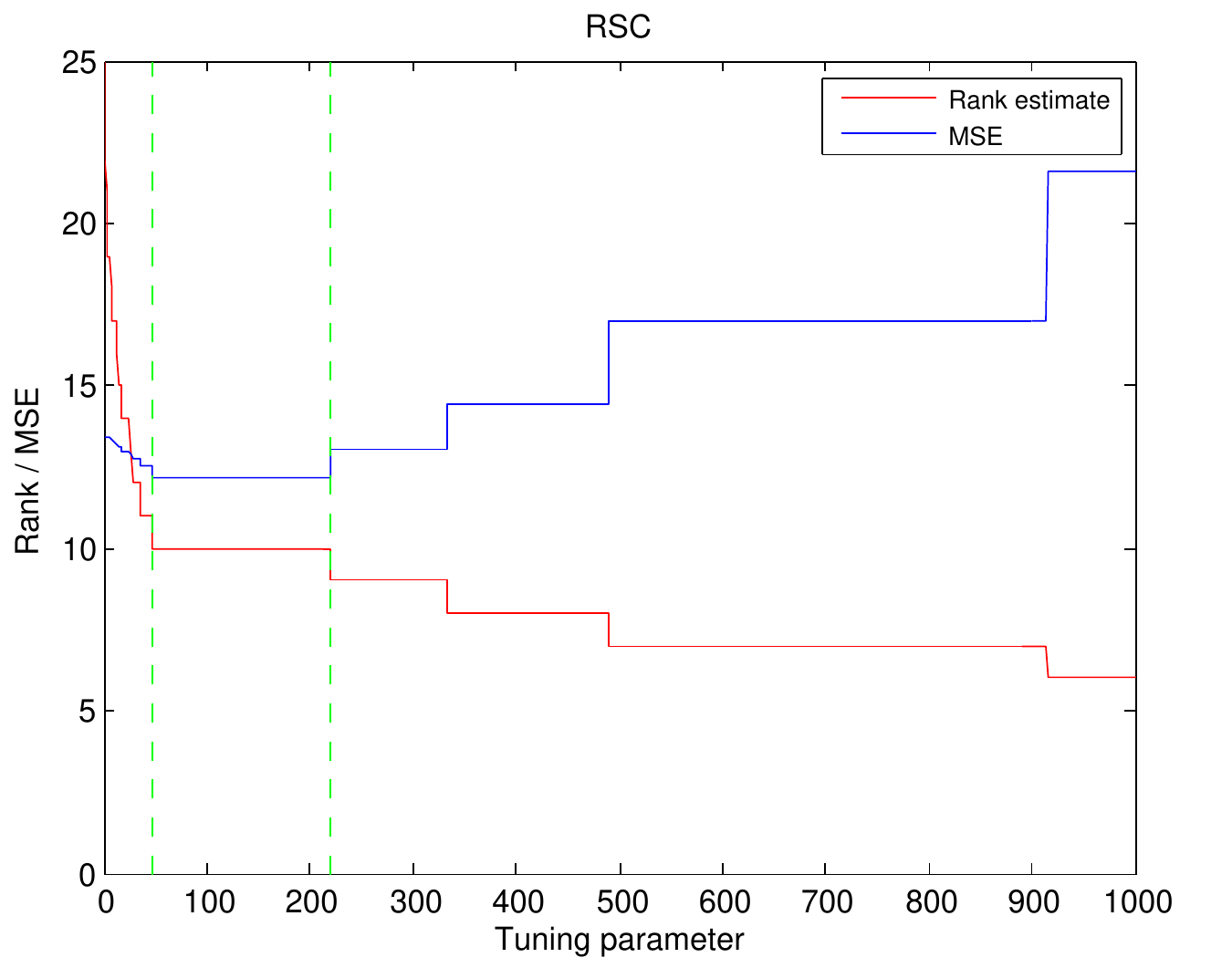}
\includegraphics[width=5cm, height=5cm ]{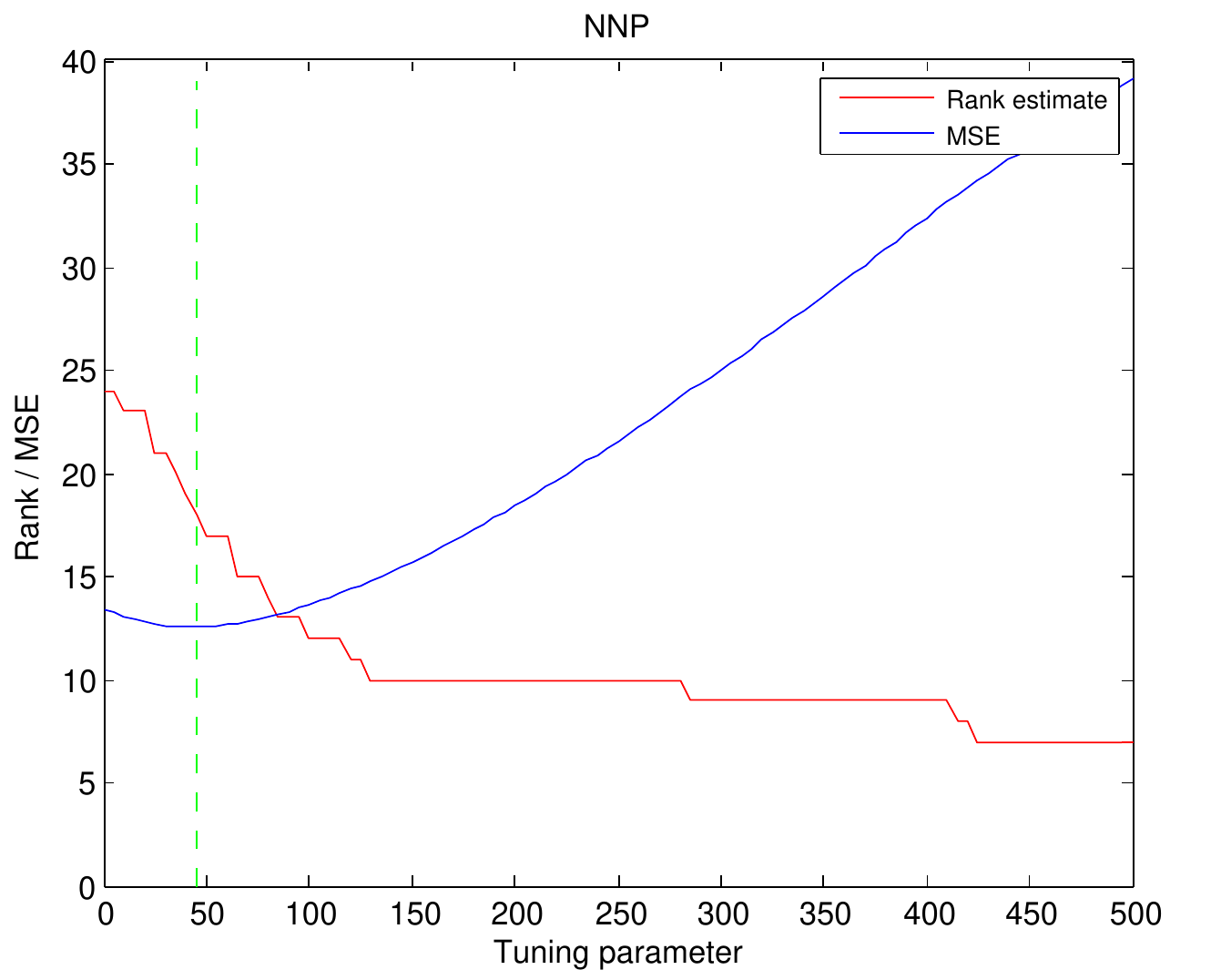}
\end{center}
\caption[]{\small{The MSE and rank of the estimators RSC (left) and NNP (right)
as a function of the tuning parameter. The rank estimate and MSE curves are plotted together for a better view of the effect of tuning on different estimation aspects.}}
\label{aap}
\end{figure}

Theorem \ref{thm:NNPoracle1} shows  that the error bounds of  the nuclear norm penalized (NNP) estimator $\widetilde{A}$ and the RSC estimator
$\wh A$ are comparable, although it is worth pointing out that our bounds for $\wh A$ are much cleaner and obtained under fewer restrictions on the design matrix.  However, there is one aspect in which the two estimators differ radically:  correct rank recovery.
 We showed in Section 2.2 that the RSC estimator corresponding to the effective value of the tuning sequence $\mu_e$ has   the correct rank and achieves the optimal  bias-variance trade-off.  This is also visible in
  the left panel of Figure \ref{aap} which shows the plots of the MSE and rank of the RSC estimate as we varied the tuning parameter of the procedure over a large grid. The numbers on the vertical axis  correspond to the range of  values of the rank  of the estimator considered in this experiment, 1 to 25. The rank of $A$ is 10.  We notice that for the same range of values of the tuning parameter, RSC has both the smallest MSE value and
the correct rank. We repeated this experiment for the NNP estimator. The right panel shows that  the smallest MSE and the correct rank are {\it not} obtained for the same value of the tuning parameter. Therefore, a different strategy for correct rank estimation via NNP is in order.
Rather than taking the rank of $\widetilde A$ as the estimator of the rank of $ A$, we consider instead, for $M = X^{\prime}X$,
\begin{eqnarray}
\label{elfje}
\tilde k= \max\{ k: \ d_k(M \widetilde A)> 2\tau  \}.\end{eqnarray}

\begin{thm}\label{elf}
Let $r=r(A)$ and assume that  $d_r(M A)> 4\tau$.
Then \[\PP\{ \tilde k\ne  r\} \le \PP\{ d_1(X' E) > \tau \} .\]
If $E$ has independent $N(0,\sigma^2)$ entries and
 $\tau= (1+\theta)  \sigma  d_1(X) (\sqrt{n}+\sqrt{q}) $,
the above  probability  is bounded by $\exp\left(- \theta^2 (n+q)/2 \right)$.
\end{thm}
\begin{proof}
After computing the sub-gradient of $f(B)=\| Y-XB\|_F^2 + 2\tau \|B\|_1$, we find that $\tilde{A}$ is a minimizer of $f(B)$ if and only if there exists a matrix $J$  with $d_1(J)\le 1$  such that
$
X^\prime X(\tilde{A} - A) = X^\prime E +\tau {U}{J} {V}^\prime $, where $\widetilde A=UDV'$ is the full SVD and $U$ and $V$ are orthonormal matrices.
The matrix ${J}$ is obtained from $D$ by setting $J_{ii}=0$ if $D_{ii}=0$ and $J_{ii}\le 1$ if $D_{ii}>0$.
Therefore,
\begin{eqnarray*}
d_1(M\tilde{A} -M A)
\le d_1(X' E) + \tau.
\end{eqnarray*}
From Horn and Johnson (1985, page 419),
\[ | d_k(M\tilde A)- d_k(MA) |\le d_1(M\tilde A-MA) \le 2\tau
\]  for all $k$, on the event $d_1(X' E)\le \tau $.
This means that $d_k(M\tilde A) > 2\tau$ for all $k\le r$ and $d_k(M\tilde A) <2\tau$ for all $k>r$, since $d_r(MA) > 4\tau$ and $d_{r+1}(MA)=0$.
The result now follows.
\end{proof}

\section{Empirical Studies}
\subsection{RSC vs. NNP}
We performed an extensive simulation study to evaluate the performance of the proposed method, RSC, and compare it with the  NNP method.
The RSC estimator  $\wh A$ was computed via the procedure outlined in Section 2.1. This method is computationally efficient  in large dimensions. Its computational complexity is the same as that of PCA.
Our choice for the tuning parameter $\mu$ was based on our theoretical findings in Section 2.
In particular, Corollary \ref{vier} and Corollary \ref{oracle1a} guarantee good rank selection and prediction performance of RSC provided that $\mu$ is just a little bit larger than $\sigma^2 (\sqrt n + \sqrt q)^2$. Under the assumption that $q<m$, we can estimate $\sigma^2$ by $S^2$; see Section \ref{subsec:unknownvar} for details. In our simulations we
 used  the adaptive tuning parameter
    $    \mu_{adap}=2 S^2 (n+q)$.
 We experimented with other constants and found that the constant equal to 2 was optimal;  constants slightly larger than 2 gave very similar results.

We compared the RSC estimator with the NNP estimator $\widetilde A$ and with the proposed
trimmed or calibrated NNP estimator, denoted in what follows by NNP$^{(c)}$.
The NNP estimator is the minimizer of the convex criterion  $\| Y - XB\|_F^2 + 2\tau \|B\|_1.$  By  the equivalent SDP characterization of the NNP-norm given in Fazel (2002),   the original minimization problem is equivalent to the convex optimization problem
\begin{eqnarray}\label{probl}
\min_{B \in {{\RR}}^{p\times n}, W_1 \in {{S}}^{n-1}, W_2 \in {{S}}^{p-1}} \| Y-XB\|_F^2+\tau ({Tr(W_1)+Tr(W_2)})\\
\hspace{-2cm} \mbox{ s.t. } \left[\begin{array}{cc} W_1 & B^T\\ B& W_2 \end{array}\right] \succeq 0.\nonumber
\end{eqnarray}
Therefore,  the NNP estimator can be computed by  adapting  the general convex optimization algorithm  {\tt SDPT3} (Toh et al. 1999) to (\ref{probl}). Alternatively, Bregman iterative algorithms can be developed; see Ma et al. (2009) for a detailed description  of the main idea. Their code, however, is specifically designed for matrix completion and does not cover the multivariate regression problem. We implemented this algorithm for the simulation study presented below.
The NNP$^{(c)}$ is our calibration of the  NNP estimator,  based on Theorem \ref{elf}.
For a given value of the tuning parameter $\tau$ we calculate the NNP estimator
$\widetilde{A}$ and obtain the rank estimate $\widetilde{r}$ from (\ref{elfje}). We then calculate the calibrated NNP$^{(c)}$  estimator as the reduced rank estimator $\wh B_{\tilde r}$, with rank equal to $\widetilde{r}$, following the procedure outlined in Section 2.1.

In our simulation study we  compared the rank selection  and the estimation performances of the
RSC estimator RSC$|_{adap}$, corresponding to  $\mu_{adap}$, with the optimally tuned RSC estimator, and the optimally tuned NNP and NNP$^{(c)}$ estimators. The last  three estimators are called RSC$|_{val}$, NNP$|_{val}$ and NNP$^{(c)}|_{val}$. They correspond to  those tuning parameters $\mu_{val}$, $\tau_{val}$ and $\tau_{val}$, respectively, that gave  the best prediction accuracy, when prediction was  evaluated on a very large independent  validation set. This comparison helps us understand  the true potential of each method in an ideal situation, and allows us to draw a stable performance comparison between the proposed adaptive RSC estimator and the best possible versions  of RSC and NNP.

We considered the following \textit{large sample-size} set up and \textit{large dimensionality} set up.
\paragraph*{Experiment 1 ($m>p$)}
We constructed the matrix of dependent variables  $X=[x_1, x_2, \cdots, x_m]^{\prime}$ by generating its rows $x_i$ as  i.i.d. realizations from a multivariate normal distribution $ \mbox{MVN}(\boldsymbol{0}, \Sigma)$,  with
$\Sigma_{jk}=\rho^{|j-k|}$,  $\rho > 0$, $1\leq j,k\leq p$.  The coefficient matrix
$A = b B_0 B_1$,  with $b > 0$,  $B_0$ is a $p\times r$ matrix and $B_1$ is a $r\times n$ matrix.  All entries  in $B_0$ and $B_1$ are i.i.d. $N(0, 1)$.  Each row in $Y = [y_1, \cdots, y_m]^{\prime}$ is then generated as $y_i = x_{i}^{\prime}A+ E_i$, $1 \leq i \leq m$, with $E_i$ denoting the $i$-th row of the noise matrix $E$ which has $m \times n$ independent $N(0,1)$ entries  $E_{ij} $.
\paragraph*{Experiment 2 ($p>m(>q)$)}
The sample size in this experiment is relatively small. $X$ is generated as  $X_0 \Sigma^{1/2}$, where $\Sigma_{jk}=\rho^{|j-k|}\in {\mathbb R}^{p\times p}$,  $X_0=X_1 X_2$, $X_1 \in {\mathbb R}^{m\times q}$, $X_2 \in {\mathbb R}^{q\times p}$ and all entries of $X_1, X_2$ are i.i.d. $N(0,1)$.  The coefficient matrix and the noise matrix are generated in the same way as in Experiment 1. Since   $p>m$, this is a much more challenging setup than the one considered in Experiment 1. Note however that  $q$, the rank of $X$,  is required to be strictly less than $m$.

Each simulated model is characterized by the following control parameters: $m$ (sample size), $p$ (number of independent variables), $n$ (number of response variables), $r$ (rank of $A$), $\rho$ (design correlation), $q$ (rank of the design), and $b$ (signal strength).
 In Experiment 1, we set $m=100,\, p=25,\, n=25,\, r=10$, and varied the  correlation coefficient $\rho=0.1, 0.5, 0.9$ and signal strength $b=0.1, 0.2, 0.3, 0.4$.
All combinations of correlation and signal strength are covered in the simulations. The results are summarized in Table \ref{simutable1a}.
In Experiment 2, we set $m=20$, $p=100$, $n=25$, $q=10$, $r=5$, and varied the correlation $\rho=0.1,\, 0.5,\, 0.9$ and signal strength $b=0.1,\, 0.2,\, 0.3$. The corresponding results are reported in Table \ref{simutable2a}.  In both tables, {MSE($A$)} and {MSE($XA$)} denote the
 $40\%$ trimmed-means of $100\cdot \|A-\hat B\|_F^2  /(pn)$  and $100\cdot \|X A-X\hat B\|_F^2  /(mn)$, respectively. We also report  the median rank estimates (RE) and the successful rank recovery percentages  (RRP). \\

\begin{table}[p]
\caption{\small{Performance comparisons of \textbf{Experiment 1}, in terms of
mean squared errors (MSE($XA$), MSE($A$)), median rank estimate (RE), and  rank recovery percentage (RRP). }}

\begin{center}
  \tiny{
\begin{tabular}{|l l|| p{1.7cm} p{1.6cm} p{1.6cm} p{1.6cm} |}

\hline
 \multicolumn{2}{|c||}{}&  {RSC$|_{adap}$}
 & {RSC$|_{val}$}
                        & {NNP$|_{val}$}
                        & {NNP$^{(c)}|_{val}$} \\

\hline
\multicolumn{6}{|c|}{$b=0.1$}\\
\hline
\multirow{2}{*}{$\rho=0.9$}
    & {MSE($XA$), MSE($A$)}  & 16.6, 5.3&  16.3, 5.2&   11.5, 3.0&  16.5, 5.3 \\
    \cline{2-6}
    & {RE, RRP} & 6, 0\% & 6, 0\% &	 12, 0\% &  6, 0\%  \\

\hline
\multirow{2}{*}{$\rho=0.5$}
    & {MSE($XA$), MSE($A$)}  &18.7, 1.4&  18.1, 1.4 &  16.2, 1.1 &  18.1, 1.4 \\
    \cline{2-6}
    & {RE, RRP} & 8, 0\% & 9, 40\% &	 16.5, 0\% &  9, 35\%  \\

\hline
\multirow{2}{*}{$\rho=0.1$}
    & {MSE($XA$), MSE($A$)}  & 19.3, 1.0& 18.0, 0.9  &   16.9, 0.8 & 18.0, 0.9  \\
    \cline{2-6}
    & {RE, RRP} & 9, 0\%& 10, 75\% &	 17, 0\% &   10, 65\%   \\

\hline
\multicolumn{6}{|c|}{$b=0.2$}\\
\hline
\multirow{2}{*}{$\rho=0.9$}
    & {MSE($XA$), MSE($A$)}  &  18.4, 7.0 & 17.9, 7.1&  15.9, 5.4 &  17.9, 7.1 \\
    \cline{2-6}
    & {RE, RRP} & 8, 0\%& 9, 20\% &	 16, 0\% &  9, 15\%  \\
\hline
\multirow{2}{*}{$\rho=0.5$}
    & {MSE($XA$), MSE($A$)}  &  16.7, 1.3 & 16.7, 1.3 &  18.9, 1.5 &  16.7, 1.3 \\
    \cline{2-6}
    & {RE, RRP} & 10, 100\% &  10, 100\% &	 19, 0\% &  10, 100\%  \\
\hline
\multirow{2}{*}{$\rho=0.1$}
    & {MSE($XA$), MSE($A$)}  &  16.5, 0.9 &  16.5, 0.9 &  19.2, 1.0 &  16.5, 0.9 \\
    \cline{2-6}
    & {RE, RRP} &  10, 100\% &  10, 100\% &	 18, 0\% &  10, 100\% \\

\hline
\multicolumn{6}{|c|}{$b=0.3$}\\
\hline
\multirow{2}{*}{$\rho=0.9$}
    & {MSE($XA$), MSE($A$)}  &  17.4, 7.0 &17.3, 6.9 &   17.7, 6.7&  17.3, 7.0 \\
        \cline{2-6}
    & {RE, RRP} &  10, 65\%  &  10, 95\% &	18, 0\% &  10, 80\%  \\

\hline
\multirow{2}{*}{$\rho=0.5$}
    & {MSE($XA$), MSE($A$)}  &  16.4, 1.3 & 16.4, 1.3 &  19.8, 1.6 & 16.4, 1.3  \\
    \cline{2-6}
    & {RE, RRP} &  10, 100\% &  10, 100\% &	 19, 0\% & 10, 100 \%  \\
\hline
\multirow{2}{*}{$\rho=0.1$}
    & {MSE($XA$), MSE($A$)}  & 16.4, 0.9 & 16.4, 0.9  &  19.9, 1.1 &  16.4, 0.9  \\
    \cline{2-6}
    & {RE, RRP} & 10, 100\%  & 10, 100\% & 19, 0\% &  10, 100\%  \\
\hline
\multicolumn{6}{|c|}{$b=0.4$}\\
\hline
\multirow{2}{*}{$\rho=0.9$}
    & {MSE($XA$), MSE($A$)}  &  16.8, 6.6 &  16.8, 6.7 &   18.7, 7.4&  16.8, 6.8  \\
    \cline{2-6}
    &{RE, RRP} &   10, 100\% &   10, 100\%&   18, 0\%&	  10, 85\% \\

\hline
\multirow{2}{*}{$\rho=0.5$}
    & {MSE($XA$), MSE($A$)}  &   16.3, 1.3 &   16.3, 1.3&   20.3, 1.7&  16.3, 1.3  \\
    \cline{2-6}
    & {RE, RRP} &  10, 100\% &  10, 100\% &	 20, 0\% &  10, 100\%  \\

\hline
\multirow{2}{*}{$\rho=0.1$}
    & {MSE($XA$), MSE($A$)}  &   16.3, 0.9 &   16.3, 0.9&  20.3, 1.1&  16.3, 0.9  \\
    \cline{2-6}
    & {RE, RRP} &  10, 100\% &  10, 100\% &	 20, 0\% &  10, 100\%  \\

\hline

\end{tabular}\\
}

\end{center}
\label{simutable1a}
\end{table}

 \begin{table}[p]

\setlength{\tabcolsep}{0.8mm}
\caption{\small{Performance comparisons of \textbf{Experiment 2}, in terms of
mean squared errors (MSE($XA$), MSE($A$), median rank estimate (RE), and rank recovery percentage (RRP).}}

\begin{center}

  \tiny{

\begin{tabular}{|l l||   p{2.0cm}  p{2.0cm} p{2.0cm} p{2.0cm} |}

\hline
 \multicolumn{2}{|c||}{}& {RSC$|_{adap}$}
                        &{RSC$|_{val}$}
                        & {NNP$|_{val}$}
                        & {NNP$^{(c)}|_{val}$}  \\

\hline
\multicolumn{6}{|c|}{$b=0.1$}\\
\hline
\multirow{2}{*}{$\rho=0.9$}
    & {MSE($XA$), MSE($A$)}  &   29.4, 3.9 &   29.4, 3.9&   36.4, 3.9 &  29.4, 3.9 \\
    \cline{2-6}
    & {RE, RRP} &  5, 100\% &  5, 100\% &	 10, 0\% &  5, 100\%  \\
    \hline

\multirow{2}{*}{$\rho=0.5$}
    & {MSE($XA$), MSE($A$)}  &  29.1, 3.9 &  29.1, 3.9 &  37.2, 3.9 &  29.1, 3.9 \\
    \cline{2-6}
    & {RE, RRP} &  5, 100\% &  5, 100\% &	 10, 0\% &  5, 100\%  \\
    \hline

\multirow{2}{*}{$\rho=0.1$}
    & {MSE($XA$), MSE($A$)}  & 29.0, 3.9 & 29.0, 3.9  &   37.2, 4.0 &  29.0, 3.9 \\
    \cline{2-6}
    & {RE, RRP} &  5, 100\% &  5, 100\% &	 10, 0\% &  5, 100\%  \\
    \hline
\multicolumn{6}{|c|}{$b=0.2$}\\
\hline
\multirow{2}{*}{$\rho=0.9$}
    & {MSE($XA$), MSE($A$)}  &   28.9, 15.7 &   28.9, 15.7&  38.7, 15.7 &  28.9, 15.7\\
    \cline{2-6}
    & {RE, RRP} &  5, 100\%&  5, 100\% &	 10, 0\% &  5, 100\%  \\
\hline

\multirow{2}{*}{$\rho=0.5$}
    & {MSE($XA$), MSE($A$)}  &  28.6, 15.7 &  28.6, 15.7 &  39.0, 15.7 & 28.6, 15.7 \\
    \cline{2-6}
    & {RE, RRP} &  5, 100\% &  5, 100\% &	 10, 0\% & 5, 100\%  \\

\hline
\multirow{2}{*}{$\rho=0.1$}
    & {MSE($XA$), MSE($A$)}  &  28.7, 15.8 &  28.7, 15.8 &  38.7, 15.8 &  28.7, 15.8\\
    \cline{2-6}
    & {RE, RRP} &  5, 100\% &  5, 100\% &	 10, 0\% &  5, 100\% \\
\hline

\multicolumn{6}{|c|}{$b=0.3$}\\
\hline
\multirow{2}{*}{$\rho=0.9$}
    & {MSE($XA$), MSE($A$)}  &  28.8, 35.3 &  28.8, 35.3 &   39.2, 35.3&  28.8, 35.3 \\
    \cline{2-6}
    & {RE, RRP}   &  5, 100\% &  5, 100\% &	10, 0\% &  5, 100\%  \\
\hline

\multirow{2}{*}{$\rho=0.5$}
    & {MSE($XA$), MSE($A$)}  &  28.5, 35.4 &  28.5, 35.4 &  39.5, 35.4 &  28.5, 35.4 \\
    \cline{2-6}
    & {RE, RRP} &  5, 100\% &  5, 100\% &	 10, 0\% & 5, 100 \%  \\
\hline

\multirow{2}{*}{$\rho=0.1$}
    & {MSE($XA$), MSE($A$)}  & 28.6, 35.5 & 28.6, 35.5  &  39.3, 35.5 & 28.6, 35.5   \\
    \cline{2-6}
    & {RE, RRP} & 5, 100\% & 5, 100\% & 10, 0\% &  5, 100\%  \\
    \hline

\end{tabular}\\
}
\end{center}
\label{simutable2a}
\end{table}

{\it Summary of simulation results.}

(i) We found  that the RSC estimator corresponding to  the adaptive choice of the tuning parameter   $\mu_{adap}=2 S^2 (n+q) \label{mu_est}$ has excellent performance. It behaves as well as the RSC estimator that uses the
 parameter $\mu$  tuned on the large validation set or the RSC estimator corresponding to the theoretical $\mu=2\sigma^2 (n+q)$.

(ii)
 When the  signal-to-noise ratio SNR := ${d_r(XA)}/{(  \sqrt q + \sqrt n)}$ is moderate or high,  with values approximately 1, 1.5 and 2, corresponding to $b = 0.2, 0.3, 0.4$,  and for low to moderate correlation   between the predictors ($\rho = 0.1, 0.5$),  RSC has excellent behavior  in terms of rank selection and means squared errors. Interestingly, NNP does not have optimal behavior in this set-up: its mean squared errors are slightly higher than those of the RSC estimator. When the noise  is  very large relative to the signal strength, corresponding to $b = 0.1$ in  Table \ref{simutable1a}, or when the correlation between  some covariates is very high, $\rho = 0.9$ in Table \ref{simutable1a},  NNP may be slightly more accurate than the RSC.

 (iii)   The NNP does not recover the correct rank, when its  regularization parameter is tuned by validation.
   Both Tables \ref{simutable1a}  and \ref{simutable2a} show that the correct rank $r$ ($r=10$ in Experiment 1 and $r=5$ in Experiment 2) is overestimated  by NNP. Our trimmed estimator,  NNP$^{(c)}$,  provides  a successful improvement over NNP in this respect.  This supports Theorem \ref{elf}.\\

 In additional simulations, we found that especially for low or moderate SNRs,
 the NNP parameter tuning problem is much more challenging than the RSC parameter tuning.
  NNP cannot    accurately estimate $A$ and consistently select the rank at the same time, for the same value of the tuning parameter. This  echoes the findings presented in  Figure \ref{aap}, and is to be expected: in NNP regularization, the threshold value $\tau$ also controls the amount of shrinkage,  which should  be mild for   large samples with relatively low contamination.  This is the case for moderate SNR and moderate correlation between predictors: the tuned $\tau$ tends to be too small,  so it cannot  introduce enough sparsity.  The same continues to be true for slightly larger values of $\tau$ that compensate for high noise level and very high correlation between predictors.   In summary,  one may not be able to build an accurate \textit{and} parsimonious  model via the NNP method, without further adjustments.

Overall, RSC is recommended over the NNP estimators, especially when we suspect
that the SNR is not very low. With large validation tuning, NNP$^{(c)}$ has the same properties as RSC -- they coincide when both methods select the same rank. But in general, the rank estimation via NNP$^{(c)}$ is much more difficult to tune and much more  computationally involved than RSC.

 For data with low SNR, an immediate extension of the RSC estimator that involves a second penalty term, of ridge-type, may induce the right amount of shrinkage needed to offset the noise in the data. This conjecture will be investigated carefully in future research.

\subsection{Tightness of the rank consistency results}
\begin{figure}[htp!]
\begin{center}
\includegraphics[width=7cm, height=7cm]{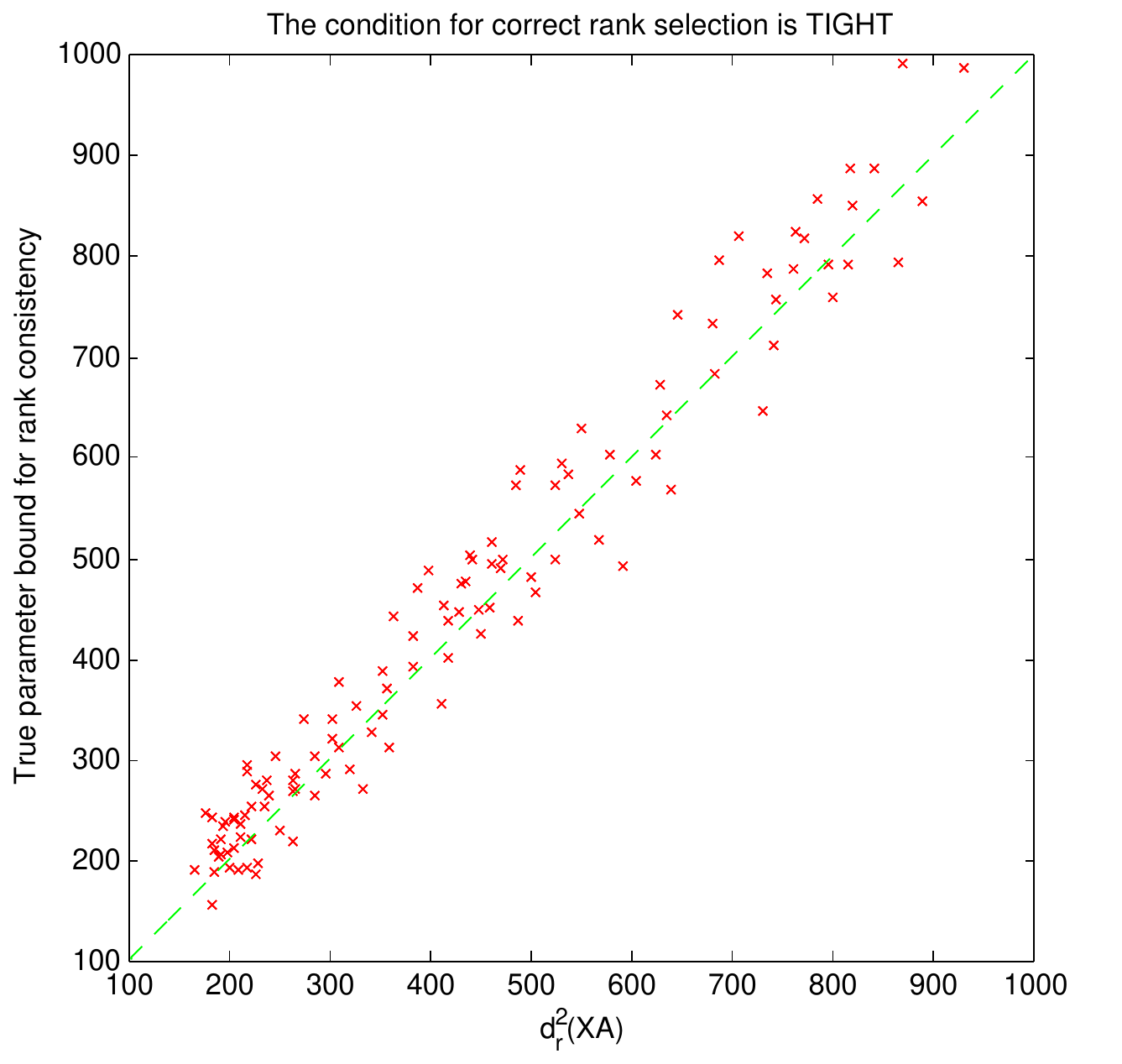}
\end{center}
\caption[]{\small{Tightness of the consistency condition.}}
\label{figconsist}
\end{figure}

It can be shown, using arguments similar to those used  in the proof of Theorem \ref{rank_regression},  that
\[ \PP\left\{ \wh k \ne r \right\} \geq
P_1 \equiv\PP\left\{ \sqrt \mu \leq d_{2r+1}(PE) \mbox{ or } d_1(PE) < \sqrt \mu-d_r(XA) \right\}.\]
On the other hand,
 the proof of Theorem \ref{rank_regression} reveals that
\[
  \PP\left\{ \wh k \ne r \right\}  \leq P_2\equiv
\PP\left\{  d_1(PE) \geq \min(\sqrt{\mu}, d_r({XA})  - \sqrt{\mu} )\right\} .
\]
Suppose now that  $2 \mu^{1/2} < d_r(XA)$ and that $r$ is small. Then $P_1$ equals $ \PP\{ d_{2r+1} (PE)\geq  \sqrt{\mu}\} $ and  is close to $P_2=\PP\{ d_1(PE) \geq \sqrt  \mu\}$ for a sparse model.
Of course, if $ \mu$ is much larger than $d_r^2(XA)$, then $P_2$ cannot be small. We use this observation to argue that, if the goal is  consistent rank estimation, then we can deviate only very little from the  requirement  $2 \mu ^{1/2}< d_r(XA)$.
This strongly suggests that the sufficient condition given in Corollary \ref{vier} for consistent rank selection is tight.
We empirically verified this conjecture for signal-to-noise ratios larger than 1 by
comparing $\mu_1=d_r^2(XA)$ with $\mu_u$, the ideal upper bound of that interval of  values of $\mu$ that give the correct rank. The value of $\mu_u$ was obtained in the simulation experiments by searching along  solution paths obtained as follows.  We constructed 100 different pairs $(X, A)$ following the simulation design outlined in the subsection above. Each pair was obtained by varying
the signal strength $b$, correlation $\rho$, the rank of $A$ and $m, n, p$.   For each  run we computed the solution path, as in Figure \ref{aap} of the previous section. From the solution path we recorded the upper value of the $\mu$ interval for which the correct rank was recovered.  We plotted  the resulting  $(\mu_1, \mu_u)$ pairs in Figure \ref{figconsist} and we conclude that the theoretical bound on $\mu$ in  Corollary \ref{vier} is  tight.

 \appendix
\section{Proof of Theorem 12} 

The starting point is the inequality
\begin{equation*}\label{part}  \| X\widetilde A - XA \|_F^2  \leq
\| XB- XA\|_F^2 +
2   \tau \left\{ \|\widetilde{A} - B\|_1 +  \| B\|_1 - \|\widetilde{A} \|_1 \right\}\end{equation*}
that holds on the event $d_1(X'E) \leq \tau$. The inequality can be deduced  from the proof of Theorem \ref{thm:NNPoracle1}.
Then,
by Lemmas 3.4  and 2.3 in Recht et al (2007) there exist two matrices $\widetilde{A}_1$ and $\widetilde{A}_2$ such that
\begin{itemize}
\item[(i)]  $ \widetilde{A} = \widetilde{A}_1 + \widetilde{A}_2$
\item[(ii)] $ r(\widetilde{A}_1) \leq 2r(B)$
\item[(iii)]  $ \| \widetilde{A} - B\|_1 =  \|\widetilde{A}_1 - B \|_1 + \|\widetilde{A}_2\|_1$
\item[(iv)]  $ \| \widetilde{A} - B\|_F^2 =  \|\widetilde{A}_1 - B \|_F^2 + \|\widetilde{A}_2\|_F^2 \ge \|\widetilde{A}_1 - B \|_F^2$
\item[(v)] $ \|\widetilde{A}\|_ 1 = \|\widetilde{A}_1\|_ 1 + \|\widetilde{A}_2\|_ 1$.
\end{itemize}
Using the display above, we find
\begin{eqnarray*}
&& \| X\widetilde A - XA \|_F^2  \\
&&\leq  \| XB-XA\|_F^2 + 2\tau \left\{ \|\widetilde{A}_1 - B\|_1 +  \|\widetilde{A}_2\|_1 +  \| B\|_1 - \|\widetilde{A}_1 \|_1  - \|\widetilde{A}_2\|_1 \right\} \\
&&\qquad  \mbox{by}  \  (i),  \ (iii) \ \mbox{and}  \ (v)\   \\
&&\leq \| XB-XA\|_F^2 + 4\tau\|\widetilde{A}_1 - B\|_1 \\
&& \leq \| XB-XA\|_F^2 + 4\tau \sqrt{r(\widetilde{A}_1 - B) }\ \|\widetilde{A}_1 - B\|_F   \mbox{ by Cauchy-Schwarz}    \\
&&\leq \| XB-XA\|_F^2 + 4\tau \sqrt{3  r(B)} \ \|\widetilde{A} - B\|_F  \mbox{ by} \ (ii)\ \mbox{and} \ (iv).
\end{eqnarray*}
Using $\lambda_{p}(M) \|\widetilde{A} - B\|_F ^2 \le \| X\widetilde{A} -X B\|_F^2$ and
$2xy\le x^2/2 + 2y^2$, we obtain
\begin{eqnarray*}
\frac{1}{2} \| X \widetilde A - XA \|_F^2  &\leq & \frac{3}{2} \| XB-XA\|_F^2 + 24 \tau^2 r(B) / \lambda_{p}(M) . \end{eqnarray*}
The proof is complete by choosing the truncated GSVD $B'$ under metric $M$, see Lemma \ref{lemma:GSVD} below.
\qed

\section{Generalized singular value decomposition}
We consider the functional
\[ G(B)= \| XB_0 - XB\|_F^2 = \text{tr}( (B-B_0)' M (B-B_0))
\] with $M=X'X=N N $ and $B_0$ is a fixed $p\times n$ matrix of rank $r$.
By the Eckhart-Young theorem, we have
the lower bound
\[ G(B) \ge \sum_{j>k} d_j^2 (XB_0) \]
for all $p\times n$ matrices $B$ of rank $k$. We now show that this infimum is achieved by the generalized singular value decomposition (GSVD) under metric $M$, limited to its $k$ largest generalized singular values.
Following Takane and Hunter (2001, pages 399-400), the GSVD of $B_0$ under metric $M$ is $UDV'$ where
$U$ is an $p\times r$ matrix, $U'MU=I_r$, $V$ is an $n\times r$ matrix, $V'V=I_r$ and $D$ is a diagonal $r\times r$ matrix, and
$ N B_0 = N UDV'. $
It can be computed via the (regular) SVD $\bar U \bar D \bar V'$ of $N B_0$. From $B_0'X' X B_0=V D^2 V'$,  the generalized singular values $d_j$ are the regular singular values of $N B_0$.
Let $B_k= U_k D_k V_k'$ by retaining as usual the first $k$ columns of $U$ and $V$.

\begin{lemma} \label{lemma:GSVD} Let $B_k$ be the GSVD of $B_0$ under metric $M$, restricted to the $k$ largest generalized singular values.
We have $$\| XB_0-XB_k\|_F^2 = \sum_{j>k} d_j^2 (X B_0). $$
\end{lemma}
\begin{proof}
Since $N B_0 = N U D V'$ and $N B_k = N U_k D_k V_k'$, we obtain
\[ \Delta= N B_0 - N B_k =N \sum_{j>k} u_j v_j' d_j =N U_{(k)} D_{(k)} V_{(k)}'\]
using the notation $U_{(k)}$ for the $p\times (r-k)$ matrix consisting of the last $r-k$ column vectors of $U$, $D_{(k)}$ is the diagonal $(r-k)\times (r-k)$ matrix based on the  last $r-k$ singular values, and
$V_{(k)}$ for the $n\times (r-k)$ matrix consisting of the last $r-k$ column vectors of $V$.
Finally,
\begin{eqnarray*}
 \| X B_0 - X B_k\|_F^2 &=& \| \Delta\|_F^2 \,=\, \| N U_{(k)} D_{(k)} V_{(k)}' \|_F^2 \\
 &=& \text{tr}\left( V_{(k)} D_{(k)} U_{(k)}' M U_{(k)} D_{(k)} V_{(k)}' \right)\\
 &=& \text{tr}\left( V_{(k) } D_{(k)} I_{(k)} D_{(k)} V_{(k)} '\right)\,
 =\, \text{tr}\left( D_{(k)}^2 \right)\,
 =\, \sum_{j>k} d_j^2 .
 \end{eqnarray*}
Recall that in the construction of the GSVD, the generalized singular values $d_j$ are the singular values of $N B_0$. Since
\[ d_j^2(N B_0) =\lambda_j ( B_0 ' M B_0 ) = \lambda_j( B_0 ' X' X B_0 ) = d_j^2 (X B_0),\]
the claim follows.
 \end{proof}

{\sc Remark.} The rank restricted estimator $\wh B_k$ given in Section 2.1 is the GSVD of  the least squares estimator $\wh B$ under the metric $M = X'X$, see Takane and Hwang (2007).

\section{Largest singular values of transformations of subGaussian matrices}

We call a random variable $W$ subGaussian with subGaussian moment $\Gamma_W $,
 if \[ \EE\left[ \exp(tW) \right] \le \exp(t^2/\Gamma_W )\] for all $t>0$. Markov's inequality implies that $W$ has Gaussian type tails:
 \[ \PP\{ |W| > t\} \le 2 \exp\{-t^2/(2\Gamma_W ) \}\]
 holds for any $t>0$. Normal $N(0,\sigma^2)$ random variables are subGaussian with $\Gamma_W=\sigma^2$.
General results on the largest singular values of matrices $E$ with subGaussian entries can be found in the survey paper by Rudelson and Vershynin (2010).
The analysis of our estimators require bounds for the largest singular values of $PE$ and $X'E$,  for which the standard results on $E$ do not apply directly.

\begin{prop}\label{een}
Let $E$ be a $m\times n$ matrix with independent subGaussian entries $E_{ij}$ with subGaussian moment $\Gamma_E $.
Let $X$ be an $m\times p$ matrix of rank $q$ and let $P=X(X'X)^{-} X'$ be the projection matrix on $R[X]$.
Then, for each $x>0$,
\begin{eqnarray*}
\PP\left\{ d_1^2( PE ) \ge 32 \Gamma_E ((n+q)\ln(5) + x) \right\} \le  2\exp \left( -x\right).
\end{eqnarray*}
In particular,
\[
\EE\left[ d_1(PE) \right] \le 15 \Gamma_E \sqrt{n+q} .\]
\end{prop}
 \begin{proof}
Let $S^{n-1}$ be the unit sphere in $\RR^n$.
First we note that
\begin{eqnarray*}
\| P E\|_2 = \sup_{u \in S^{p-1},\ v\in S^{n-1} } < Pu, Ev > = \sup_{u\in U,\ v\in S^{n-1} } < u, Ev>
\end{eqnarray*}
with $U= P S^{p-1} = \{ u= Ps:\ s\in S^{p-1} \}$.
Let $\M$ be a $\delta$-net of $U$ and $\N$ be a $\delta$-net for $S^{n-1}$ with $\delta=1/2$.   Since the dimension of $U$ is $q$ and$\|u\|\le 1$ for each $u\in U$, we need at most $5^{q}$  elements in $\M$ to cover $U$ and $5^n$ elements to cover $S^{n-1}$,
see Kolmogorov and Tikhomirov (1961). A standard discretization trick, see, for instance,  Rudelson and Vershynin (2010, proof of Proposition 2.4), gives
\[ \| PE\|_2 \le 4 \max_{ u\in \M,\ v\in \N } <  u, Ev>.\]
Next, we
write
$ < u, Ev> = \sum_{i=1}^m u_i <E_i,v> $
and note that
each $<E_i,v>$ is subGaussian with moment $\Gamma_E $, as
\[ \EE\left[ \exp(t<E_i,v> ) \right] = \prod_{j=1}^n  \EE\left[ \exp(t v_j E_{ij} ) \right]  \le \exp(t^2 \sum_j v_j^2 / \Gamma_E  ) = \exp(t^2/\Gamma_E  ).
\]
It follows that each term in  $\sum_{i=1}^m  u_i <E_i,v> $ is subGaussian, and
$< u, Ev> $ is subGaussian with subGaussian moment $ \Gamma_E  \sum_{i=1}^m u_i^2 = \Gamma_E$. This implies the tail bound
\[ \PP\{ | < u, Ev > | > t \} \le 2 \exp\{ -t^2 / (2\Gamma_E ) \}\]
for each fixed $u$ and $v$ and all $t>0$.
Combining the previous two steps, we obtain
\[ \PP\left\{ \| PE\|_2 \ge 4t \right\} \le 5^{ n+ q} 2\exp\{ -t^2 / (2\Gamma_E   ) \}\]
for all $t>0$. Taking $t^2 =  2\{ \ln(5) (n+ q) + x \} \Gamma_E $ we obtain the first claim.
The second claim follows from  this tail bound.
 \end{proof}

 \section{Auxiliary results}

\begin{lemma}\label{basis1}
Let $X$ be a non-negative random variable with $\EE[X]=\mu$ and
$ \PP\{ X-\mu \ge t \} \le \exp(-t^2/2) $
for all $t\ge0$. Then
we have
\[ \EE\left[ X^2   \right] \le \mu^2 + \mu \sqrt{2\pi} +2.
\]
Moreover,
 for any $\xi>0$, we have
\[ \EE\left[ \left (X^2 - (1+\xi)^2 \mu^2 \right)_{+} \right] \le 2 (1+{\xi}^{-1}) \exp(-\xi^2\mu^2/2).
\]
\end{lemma}
\begin{proof}
The following string of inequalities are self-evident:
\begin{eqnarray*}
\EE\left[ X^2  \right] &=& \int_0^\infty \PP\{ X^2 \ge x\}\, {\rm d}x
\, \le\, \mu^2 + \int_{\mu}^\infty 2x \PP\{ X \ge  {x} \} \, {\rm d}x\\
&\le& \mu^2 + \int_0^\infty  2(x+\mu) \exp\left( -\frac12 x^2 \right) \, {\rm d}x\,
=\, \mu^2 + \mu \sqrt{2\pi} +2.
\end{eqnarray*}
This proves our first claim.
The second  claim is easily deduced as follows:
\begin{eqnarray*}
\EE \left[ \left (X^2 - (1+\xi)^2 \mu^2 \right)_{+} \right]
&\le& \EE \left[X^2   1_{\left \{ X\ge (1+\xi)\mu\right\} } \right] \,
=\, \int_{(1+\xi)\mu} ^\infty 2t\PP\{ X\ge t\}\, {\rm d}t\\
&\le& (1+\xi^{-1})  \int_{\xi\mu} ^\infty 2t\exp(-t^2/2) \, {\rm d}t\\
&=& 2  (1+\xi^{-1}) \exp(-\xi^2\mu^2/2).
\end{eqnarray*}
The proof of the lemma is complete.
\end{proof}

\begin{lemma} \label{chisq}
Let $Z_d$ be a $\chi^2_d$ random variable with $d$ degrees of freedom. Then
\begin{eqnarray*}
\PP\left\{ Z_d-d \le - x\sqrt{2d} \right\} \le \exp\left(- \frac{x^2}{ 2+2x\sqrt{2/d}}. \right)
\end{eqnarray*}
In particular, for any $0<t<1$,
\begin{eqnarray*}
\PP\left\{ Z_d \le (1-t)d \right\} \le \exp\left\{- t^2 d / 4(1+t)\right\}.
\end{eqnarray*}
\end{lemma}
\begin{proof}
See Cavalier et al (2002, page 857) for the first claim. The second claim follows by taking $x=t(d/2)^{1/2}$.
\end{proof}

{\bf Acknowledgement.} We would like to thank Emmanuel Cand\`es, Angelika Rohde and Sasha Tsybakov for stimulating conversations in Oberwolfach, Tallahassee and Paris, respectively. We also thank the associate editor and the referees for their constructive remarks.

\end{document}